\documentclass[[11pt,a4paper]{amsart}
\usepackage{amssymb}
\usepackage{color}
\usepackage[arrow, matrix, curve]{xy}
\usepackage{hyperref}

\newtheorem{thm}{Theorem}
\newtheorem{lem}{Lemma}
\newtheorem{cor} {Corollary}

\newtheorem{pro} {Proposition}
\newtheorem{df} {Definition}
\newtheorem{obs}{Observation}
\newtheorem{rmk}{Remark}
\newcommand {\sima}{\parallel M_1,\partial M_1\parallel}
\newcommand {\simb}{\parallel M_2,\partial M_2\parallel}
\newcommand {\glu}{M_1\cup_f M_2}

\newcommand {\m} {C_*\left(M,M^\prime\right)\otimes_{{\Bbb Z}G}{\Bbb Z}}
\newcommand {\kk} {C_*\left(K,K^\prime\right)\otimes_{{\Bbb Z}G}{\Bbb Z}}
\newcommand {\lt} {C_*\left(L,L^\prime\right)\otimes_{{\Bbb Z}G}{\Bbb Z}}
\newenvironment{pf}{{\it Proof:}\quad}{\hfill QED}
 
\begin{document}
 
\title{Multicomplexes, bounded cohomology and \\
additivity of simplicial volume }
\author{Thilo Kuessner}
\date{}
\maketitle
\begin{abstract}
We discuss some additivity properties of the simplicial volume for manifolds with boundary: we give
proofs of
additivity for glueing amenable boundary components and of superadditivity for glueing
amenable submanifolds of the boundary,
and we discuss doubling of 3-manifolds.
\end{abstract}
 
\noindent                                                
This paper is devoted to the behaviour of simplicial volume under cutting and
pasting along amenable submanifolds. Such results have been proved by Gromov in \cite{gro}, 
for closed manifolds (\cite{gro}, Section 3) by fairly elementary arguments,
and for open 
manifolds (\cite{gro}, Section 4) by much more advanced arguments. The aim of this article 
is to write complete proofs for the case of compact manifolds 
with boundary, using elementary arguments which are
closer to the arguments that were needed for the closed case.

The simplicial volume is a homotopy invariant of compact manifolds. It was defined in \cite{gro}, using the $l^1$-norm on singular chains, as follows. \\
\\
{\bf Definition:} 
{\em a) Let $\left(X,Y\right)$ be a pair of spaces, and $h\in
H_*\left(X,Y;{\Bbb R}\right)$ a relative homology class. Its Gromov norm $\parallel h\parallel$ is defined by}
$$\parallel h\parallel:=
\inf\left\{\sum_{i=1}^r\mid a_i\mid :
\sum_{i=1}^r a_i\sigma_i\mbox { represents }h\right\}.$$
{\em b) Let $M$ be 
a compact, connected, orientable $n$-manifold, Let $\left[M,\partial M\right]\in H_n\left(M,\partial M;{\Bbb R}\right)$ 
be a fundamental class, that is the image of a generator of $H_n\left(M,\partial M;{\Bbb Z}\right)$ under the canonical 
homomorphism $H_n\left(M,\partial M;{\Bbb Z}\right)\rightarrow H_n\left(M,\partial M;{\Bbb R}\right)$. 
Define the simplicial volume of $M$ to be the Gromov norm of the fundamental class:}
$$\parallel M,\partial M\parallel:=\parallel\left[M,\partial M\right]\parallel.$$

(This definition extends to non-connected manifolds by summing the simplicial volumina 
of the connected
components, and to non-orientable manifolds by taking half of the simplicial volume of the oriented double covering.)
 
The simplicial volume remains an invariant about which relatively little is known. This article is devoted to the study of the behaviour of simplicial volume with
respect to cut and paste.
That means, we are given a properly embedded $(n-1)$-submanifold $F\subset M$ (that is $F\cap\partial M=\partial F$), and we wish to compare $$\parallel M_F,\partial M_F\parallel
\mbox{\ to\ } \parallel M,\partial M\parallel.$$

Here, $M_F$ denotes the manifold obtained by cutting off $F$, that is $$M_F:=\overline{M - N\left(F\right)}$$ for a regular neighborhood $N\left(F\right)$ of $F$.
 
It may happen that $\parallel M,\partial M\parallel < \parallel M_F,\partial M_F\parallel$,
as there may be fundamental cycles of $M$ which are not the images of fundamental cycles of $M_F$. (This is because there are fundamental cycles of $M$ using simplices that intersect $F$ transversally.) For example, we showed in \cite{k} that $$\parallel M,\partial M\parallel < \parallel M_F,\partial M_F\parallel$$
if int$\left(M\right)$ is a hyperbolic $n$-manifold of finite volume, $n\ge3$, 
and $F$ a closed geodesic hypersurface.                      

Somewhat counter-intuitively, there are also examples satisfying a strict inequality $$\parallel M,\partial M\parallel
 > \parallel M_F,\partial M_F\parallel,$$ which is possible because fundamental cycles for $M_F$ need not fit
together at the copies of $F$.

We are going to consider the case that $\pi_1 F$ is amenable. Slightly more generally we will prove:

\begin{thm}\label{Theorem1}Let $M_1,M_2$ be two compact $n$-manifolds,
$$A_1\subset \partial M_1, A_2\subset\partial M_2$$
properly embedded $(n-1)$-dimensional submanifolds,
$f:A_1\rightarrow A_2$ a homeomorphism, and
$$M=\glu$$ the glued manifold. Assume that $\pi_1\partial M\rightarrow\pi_1M$ is injective.
If $$im\left(\pi_1A_1\rightarrow\pi_1M_1\right), im\left(\pi_1A_2\rightarrow\pi_1M_2\right)$$ are amenable and $f_*$ restricts to an isomorphism
$$f_*:ker\left(\pi_1A_1\rightarrow\pi_1M_1\right)\rightarrow
ker\left(\pi_1A_2\rightarrow\pi_1M_2\right), ker\left(\pi_1\partial A_1\rightarrow\pi_1M_1\right)\rightarrow
ker\left(\pi_1\partial A_2\rightarrow\pi_1M_2\right)$$
then $$\parallel M,\partial M\parallel\ge\sima+\simb.$$
If moreover $A_1,A_2$ are connected components of $\partial M_1$ resp.\ $\partial
M_2$,
then $$\parallel M,\partial M\parallel=\sima+\simb.$$\end{thm}

We prove analogous facts if $A_1,A_2$ are in the boundary of the same manifold $M_1$.\\        

Bucher et al. (\cite[Theorem 1.3]{bbfipp}) use other methods from \cite{iva} to give an easier proof for the equality $\parallel M,\partial M\parallel=\sima+\simb$ under the slightly stronger assumption that all components of $\partial M$ have amenable fundamental group. In our result, one only needs the boundary components which are actually glued to be amenable, and also one only needs amenability of $im(\pi_1A_i\to\pi_1M_i)$. We remark that in view of the counterexamples in \hyperref[examples]{Subsection \ref*{examples}} neither the assumption $ker\left(\pi_1A_1\rightarrow\pi_1M_1\right)\simeq
ker\left(\pi_1A_2\rightarrow\pi_1M_2\right)$ nor the assumption of the $A_i$ being connected can be avoided.

Applied to 3-manifolds $M$, Theorem 1 means that simplicial volume is additive with respect to glueing along incompressible tori and superadditive with respect to glueing along incompressible annuli. If $\partial M$ consists of tori,  Soma proved in
\cite{som1} that simplicial volume is
additive with respect to glueing incompressible tori or annuli. We think that our proof, apart from being a generalisation to manifolds with arbitrary boundary,
should be of interest because the proof in \cite{som1} heavily relies on Theorem 6.5.5.\ from Thurston's lecture notes, of which no published proof is available so far.\\

To give some non-rigorous motivation for the arguments,
let us consider a toy example: the wedge $M_1\vee M_2$ of two closed manifolds $M_1$ and $M_2$.
(This is not a manifold but one may define 
a fundamental class in the obvious way.) Here
$\parallel M_1\vee M_2 \parallel \le \parallel M_1\parallel +\parallel M_2 \parallel$ trivially holds. We would like to show 
the opposite inequality
$\parallel M_1\vee M_2\parallel\ge\parallel M_1\parallel +\parallel
M_2\parallel$. That means that we have to find an efficient way to map representatives of $\left[M_1\vee M_2\right]$
to a sum of representatives of $\left[M_1\right]$ and of $\left[M_2\right]$.
So
we look
for a
chain map $r$, leftinverse to the inclusions, which maps
simplices in $M_1\vee M_2$
to simplices either in $M_1$ or in $M_2$.
 
The universal cover of $M_1\vee M_2$ is a tree-like complex made from copies of $\widetilde{M_1}$ and $\widetilde{M_2}$. In a tree, any 
triple of vertices in general position has a unique central point, belonging to all 
the three (unique) geodesics connecting the vertices.
More generally, in a tree-like complex, one might try to construct a 'central' simplex associated to any 'nondegenerate' tuple of at least three points.
For example, if $M_1$ and $M_2$ admitted hyperbolic metrics, one would have unique geodesics connecting
any two given
vertices in $\widetilde{M_1}$ or $\widetilde{M_2}$.
Thus there are also unique geodesics
in the tree-like complex, and one can actually show that, for any
'nondegenerate' tuple of vertices, the associated set of geodesics intersects the full 1-skeleton of exactly one top-dimensional simplex. This 'central' 
simplex belongs to a copy of either $\widetilde{M_1}$ or $\widetilde{M_2}$. The chain map $r$ is defined by mapping each simplex $\sigma$ in $M_1\vee M_2$ to the
projection of the 'central' simplex which is associated
to the vertices of some lift $\tilde{\sigma}$ of $\sigma$ (resp.\ to zero if 
the tuple of vertices is 'degenerate', i.e.\ if no 'central' simplex
exists).
 
There is clearly no such construction for arbitrary manifolds. However, using bounded cohomology 
one can reduce the proof of Theorem 1 to aspherical, minimally complete multicomplexes, and we show in 
Chapter 2 that 'central' simplices can be constructed if $M_1$ and $M_2$ were
aspherical,
minimally complete multicomplexes. (Such multicomplexes have in fact several features in common with hyperbolic spaces.)

What we actually prove in Chapter 2 is a construction of 'central' simplices in more general tree-like complexes,
namely in
multicomplexes glued
along a submulticomplex $A$, with a 
group $G=\Pi\left(A\right)$ acting with certain transitivity properties on $A$.

It might be tempting to consider the quotients with respect to the $G$-action to reduce 
the glueing to a wedge. However, this would raise several technical problems. First, these quotients 
are not multicomplexes, since an n-simplex need not have 
n+1 distinct vertices. Second, more seriously, these quotients are not minimal: consider, for example, two 1-simplices which are homotopic (rel.\ boundary)
to concatenations $k*a_1*l$ resp.\ $
k*a_2*l$, where $a_1,a_2$ are 1-simplices in $A$ and $k,l$ are 1-simplices not in $A$. Then these 1-smplices are homotopic 
(rel.\ boundary) but distinct in $K/\Pi\left(A\right)$. Third, probably the most serious problem, the fundamental group of the quotient complex is $\pi_1K/N\left(\pi_1A\right)$, where $N\left(\pi_1A\right)$ is the normalizer of $\pi_1A$ and need not be amenable even if $\pi_1A$ is amenable. But if the deck group is not amenable, the projection need not be an 
isometry of Gromov norms. For these reasons we work with quotients of chain complexes (in the form of tensor products) rather with quotients of multicomplexes.\\

There are two approaches to bounded cohomology, namely Gromov's approach via multicomplexes that is used in this paper and Ivanov's approach via strong resolutions in \cite{iva}. In the meantime since the first version of this paper was written, the authors of \cite{bbfipp} have worked out an approach to the questions discussed in this paper via Ivanov's method. While their assumption is slightly more restrictive than ours (they have to assume that all boundary components are amenable, not just the glued ones), that more abstract method certainly allows a more concise approach to this kind of questions. Also some of the results of Gromov's discussed in the introductory Section 1 of this paper have more concise proofs via Ivanov's approach, cf.\ the discussion in the proof of Proposition 3. Still we think that it is helpful to have an alternative approach, which geometrically might be more intuitive than the one using abstract resolutions.

Convention: 
For simplicity, we assume all manifolds to be oriented. All theorems generalise in an obvious way to non-orientable manifolds.

I would like to thank Michel Boileau and Elmar Vogt for discussions about the paper.

\section{Multicomplexes}
In this introductory section we will collect facts about multicomplexes and bounded cohomology from Gromov's work \cite{gro} which we will need lateron. One should be aware that the results in this section about bounded cohomology have by now more concise (and complete) proofs by the work of Ivanov \cite{iva} (see also \cite{bbfipp}) but we will present here Gromov's approach with the aim to introduce his machinery of multicomplexes which we are going to use in our arguments lateron.
\subsection{Definitions}
 
\begin{df} 
Let $V$ be a set. A {\em multicomplex} $K=(V,S,\left\{\partial_i\right\}_{i\in\bf N})$ with vertices $V$ is given by the following data:
\begin{itemize}
\item[-] a (possibly empty) set $S_{v_0,\ldots,v_n}$ for each $n\in\bf N$ and each ordered subset $(v_0,\ldots,v_n)\in V^{n+1}$ (the "$n$--simplices with vertices $v_0,\ldots,v_n$"),\\
\item[-] a map $\partial_i:S_{v_0,\ldots,v_n}\rightarrow S_{v_0,\ldots,\hat{v}_i,\ldots,v_n}$ for each $n\in\bf N$, each $i\in\left\{0,1,\ldots,n\right\}$ and each $(v_0,\ldots,v_n)\in V^{n+1}$.
\end{itemize}
%
\end{df}

A {\em simplicial map} from a multicomplex $K=(V_K,S,\left\{\partial_i\right\})$ to a multicomplex $K^\prime=(V_{K^\prime},S^\prime,{\{\partial^\prime_i\}})$ is given by a map $f:V_K\rightarrow V_{K^\prime}$ and a map $F:S_{v_0,\ldots,v_n}\rightarrow S^\prime_{f(v_0),\ldots,f(v_n)}$ for each $(v_0,\ldots,v_n)\in V_K^{n+1}, n\in\bf N$, such that $\partial_i^\prime F=F\partial_i$ for all $i\in\left\{0,1,\ldots,n\right\}$.

A {\em subcomplex} $K^\prime\subset K$ of a multicomplex $K=(V_K,S,\left\{\partial_i\right\})$ is a multicomplex $K^\prime=(V_{K^\prime},S^\prime,{\{\partial^\prime_i\}})$ with $V_{K^{\prime}}\subset V_K$, $S^\prime_{v_0,\ldots,v_n}\subset S_{v_0,\ldots,v_n}$ for all $(v_0,\ldots,v_n)\in V_{K^\prime}^{n+1}$, $n\in\bf N$, and $\partial_i^\prime:S^\prime_{v_0,\ldots,v_n}\rightarrow S^\prime_{v_0,\ldots,\hat{v}_i,\ldots,v_n}$ the restriction of $\partial_i:S_{v_0,\ldots,v_n}\rightarrow S_{v_0,\ldots,\hat{v}_i,\ldots,v_n}$ to $S^\prime_{v_0,\ldots,v_n}$ for all $(v_0,\ldots,v_n)\in V_{K^\prime}^{n+1}$, $n\in\bf N$, $i\in\left\{0,1,\ldots,n\right\}$. $\left(K,L\right)$ is a {\em pair of multicomplexes} if $K$ is a multicomplex and $L$ is a subcomplex of $K$.

The geometric realisation $\mid K\mid$ of a multicomplex $K=(V,S,{\{\partial_i\}})$ is $$\mid K\mid=\bigcup_{n\in \bf N}\bigcup_{(v_0,\ldots,v_n)\in V^{n+1}}S_{v_0,\ldots,v_n}\times\Delta^n/\sim,$$
where $\Delta^n$ denotes the standard $n$--simplex and the equivalence relation
$\sim$ is
generated by the relations $$(\partial_i\sigma,x)\sim(\sigma,j_i(x))$$
for all $n\in{\bf N},(v_0,\ldots,v_n)\in V^{n+1},\sigma\in S_{v_0,\ldots,v_n},i\in\left\{0,\ldots,n\right\},x\in\Delta^{n-1}$, where $j_i:\Delta^{n-1}\rightarrow\Delta^n$ denotes the standard inclusion as $i$-th face. The topology on $\mid K\mid$ is determined by the 
condition that a set $A\subset\mid K\mid$ is closed 
if and only if its intersection with each simplex is closed. 

\begin{df} A multicomplex $K$ is minimally complete, 
if the following holds:
whenever $\sigma:\Delta^n\rightarrow \mid K\mid$ is a singular n-simplex, whose restriction to the boundary $f |_{\partial\Delta^n}:\partial\Delta^n\rightarrow\mid K\mid$ is a simplicial embedding, then $\sigma$ is homotopic relative $\partial\Delta^n$ to
a {\em unique} simplex in $K$.\end{df}
\noindent
For a multicomplex $K$, its {\em $n$-skeleton} $K_n$ is the union of its simplices of dimensions $\le n$. For a simplex $\sigma\in K$, its {\em $n$-skeleton} $\sigma_n$ is the union of its iterated faces of dimensions $\le n$.

\begin{df}\label{aspherical} A minimally complete
multicomplex $K$ is aspherical if $\sigma_1\not=\tau_1$ for all $\sigma\not =\tau\in K$.\end{df} 

{\bf Orientations.} Let $\sigma:\Delta^n\rightarrow \mid K\mid$ be a simplex of $K$. Each permutation $\pi$ of the vertices extends to a unique affine homeomorphism $f_\pi:\Delta^n\rightarrow \Delta^n$. We say that $\sigma$ and $\sigma\circ f_\pi$ have the same orientation if
$\pi$ is the product of an even number of transpositions, and the opposite orientation else.

Let $\pi=\left(01\right)$ be the transposition of the $0$-th and $1$-th vertex of $\Delta^n$. For an $n$-simplex $\sigma$,
we will denote $\sigma\circ f_{\left(0,1\right)}$ by $\overline{\sigma}$ and we will say that {\em $\overline{\sigma}$ is the simplex $\sigma$ with the opposite orientation}.

{\bf Group Actions.} A group $G$ {\em acts simplicially} on a pair of multicomplexes $\left(K,L\right)$ if it acts on the set of simplices of $K$, mapping simplices in $L$ to simplices in $L$, such that the action commutes with all face maps, that is $g\partial_i\sigma=\partial_i\left(g\sigma\right)$ for all $g\in G, \sigma\in K, i\in\left\{0,\ldots,dim\left(\sigma\right)\right\}$. For $g\in G$ and $\sigma$ a simplex in $K$, we denote by
$g\sigma$ the simplex obtained by this action.

For a connected multicomplex $K$ we write
$\pi_1K$ for $\pi_1\left(\mid K\mid,p\right)$ with some $p\in\mid K\mid$. 

{\bf Simplicial Approximation.} The Simplicial Approximation Theorem holds for maps between multicomplexes. Indeed, the proof of the Simplicial Approximation Theorem for maps between simplicial complexes makes no use of the assumption that there are no $n$-dimensional simplices with the same $n-1$-skeleton.

\subsection{Bounded Cohomology}

For a multicomplex $K$, let $F_j\left(K\right)$ be the ${\Bbb R}$-vector space with basis 
the set of $j$-simplices of $K$. Let $O_j\left(K\right)$ be the subspace generated by the set 
$\left\{\sigma+sign\left(\pi\right) \sigma\circ f_\pi\right\}$ where $\sigma$ runs over all $j$-simplices, $\pi$ runs over all permutations of the vertices of $\Delta^n$.

Define the {\em $j$-th chain group} $C_j\left(K\right)=F_j\left(K\right)/O_j\left(K\right).$
Note that for any $j$-simplex $\sigma$ in $K$ we have $\overline{\sigma}=-\sigma$ in $C_j\left(K\right).$

Define the {\em group of antisymmmetric cochains} $C^j_a\left(K\right):=Hom_{\Bbb R}\left(C_j\left(K\right),{\Bbb R}\right).$

If $\left(K,L\right)$ is a pair of multicomplexes, we get an inclusion $C_j\left(L\right)\rightarrow
C_j\left(K\right)$ and we define 
$$C^j_a\left(K,L\right):=\left\{\omega\in C^j_a\left(K\right): \omega\left(c\right)=0 \mbox{ for all } c\in C_j\left(L\right)\right\}$$ with the norm 
$$\parallel \omega\parallel_\infty:=\sup\left\{\omega\left(\left[\sigma\right]
\right) : \sigma \mbox{ j-simplex}\right\}.$$
Let
$$C_b^j\left(K,L\right):=\left\{\omega\in C^j_a\left(K,L\right):
\parallel\omega\parallel_\infty < \infty\right\}.$$
The coboundary
operator 
preserves $C_b^*\left(K,L\right)$, hence induces
maps $\delta_b^j:C_b^j\left(K,L\right)\rightarrow
C_b^{j+1}\left(K,L\right)$. Define the {\em bounded cohomology}
of $\left(K,L\right)$ by $H_b^j\left(K,L\right):=
ker\delta_b^j/im\delta_b^{j-1}$. The norm $\parallel .\parallel_\infty$
induces a pseudo-norm $\parallel .\parallel$ on
$H_b^*\left(K,L\right)$.

Before introducing the bounded cohomology of spaces, let us make the following remark.
\begin{rmk}
Let $\left(X,Y\right)$ be a pair of topological spaces such that each path-component of $X$ resp.\ $Y$ contains infinitely many points. Let $S_*^{sing}\left(X\right)$ resp.\ $S_*^{sing}\left(Y\right)$ be the complexes of singular simplices. Let $\hat{S}_*^{sing}\left(X\right)\subset S_*^{sing}\left(X\right)$ resp.\ $\hat{S}_*^{sing}\left(Y\right)\subset S_*^{sing}\left(Y\right)$ be the subcomplexes consisting of
those {\em singular simplices 
which have pairwise distinct vertices}. Then the inclusion $$incl:\left(\hat{S}_*^{sing}\left(X\right),\hat{S}_*^{sing}\left(Y\right)\right)\rightarrow \left(S_*^{sing}
\left(X\right),S_*^{sing}\left(Y\right)\right)$$ 
induces an isometric isomorphism in homology, thus 
we can compute the Gromov norm of homology classes by considering cycles in 
$\left(\hat{S}_*^{sing}\left(X\right),\hat{S}_*^{sing}\left(Y\right)\right)$ only.\end{rmk}
\begin{pf}To prove surjectivity of $incl_*$ in homology let $\sum_{i=1}^ra_i\sigma_i$ 
be a (relative) cycle representing a given homology class $h$. Then 
form a simplicial complex $P$ by taking disjoint copies of the standard simplex $\Delta_1,\ldots,\Delta_r$ and 
identifying $\partial_i\Delta_k$ to $\partial_j\Delta_l$ if and only if $\partial_i\sigma_k=\partial_j\sigma_l$. 
Clearly there is a continuous map $f\colon (P,\partial P)\to (X,Y)$ realising $h$ in the 
sense that $f_*\left[k\right]=h$ for the class $k\in H_*(P,\partial P)$ represented by $\sum_{i=1}^ra_i \Delta_i$. (If $h$ is an integral homology class, then $P$ is a triangulated 
pseudomanifold, see e.g.\ \cite[Section 2]{bfp}. For a real homology class, in general $P$ will just be a simplicial complex, not a pseudomanifold.) 
 
Since 
each component of $X$ and $Y$ has infinitely many points, we can homotope 
$f\mid_{P_0}$ (the restriction of $f$ to the $0$-skeleton) to an injective map $g\colon (P_0,(\partial P)_0)\to (X,Y)$. 

Since the inclusion $P_0\to P$ of the $0$-skeleton $P_0$ in $P$ is a 
cofibration, this homotopy extends to $P$ and we get a map $g\colon (P,\partial P)\to 
(X,Y)$ which represents an element $\sum_{i=1}^ra_i\sigma_i^\prime\in \left(\hat{S}_*^{sing}\left(X\right),\hat{S}_*^{sing}\left(Y\right)\right)$ that is (relatively) homotopic and hence (relatively) homologous to $\sum_{i=1}^ra_i\sigma_i$ and thus represents $h$. This shows surjectivity 
of $H_*\left(\hat{S}_*^{sing}\left(X\right),\hat{S}_*^{sing}\left(Y\right);{\mathbb R}\right)\to H_*\left({S}_*^{sing}\left(X\right),{S}_*^{sing}\left(Y\right);{\mathbb R}\right)$. 

Similarly, if a (relative) cycle $z\in \left(\hat{S}_*^{sing}\left(X\right),\hat{S}_*^{sing}\left(Y\right)\right)$ is the (relative) boundary of a chain $w\in \left({S}_*^{sing}\left(X\right),{S}_*^{sing}\left(Y\right)\right)$, then the same argument yields a homotopy 
(rel.\ $\partial w$) from $w$ to a chain with pairwise distinct vertices and still bounded by $z$, which proves injectivity of $incl$. Obviously, the construction does not increase norms. 

See also \cite[Lemma 2a]{k}.\end{pf}\\


We then
define
the {\em bounded cohomology of (pairs of) topological spaces} as follows: for a topological space $X$ consider 
the multicomplex $S_*\left(X\right):=\hat{S}_*^{sing}\left(X\right)$, which is formed by all singular simplices with pairwise distinct vertices
\footnotemark\footnotetext[1]{The topology of this multicomplex is given by the condition that a set is closed if and only if its intersection with each simplex is closed.}, and the usual face maps $\partial_i$. 
Then we let $C^j\left(X\right)
:=Hom_{\Bbb R}\left(F_j\left(S_*\left(X\right)\right),{\Bbb R}\right)$ 
for each $j\in{\Bbb N}$. If $Y$ is a subspace of $X$, we define $$C^j\left(X,Y\right)
:=\left\{\omega\in C^j\left(X\right): \omega\left(c\right)=0 
\mbox{\ for all\ }c\in C_j\left(Y\right)\right\}
$$ with the norm $$\parallel \omega\parallel_\infty:=\sup\left\{\omega\left(\sigma\right): \sigma\mbox{\ simplex in\ }S_*\left(X\right)\right\}$$ and we let $$
C_b^j\left(X,Y\right):=\left\{\omega\in C^j\left(X,Y\right):\parallel
\omega\parallel_\infty < \infty\right\}.$$
Again, we have the coboundary operator
$\delta_b^j:C_b^j\left(X,Y\right)\rightarrow C_b^{j+1}\left(X,Y\right)$
and we let $H_b^j\left(X,Y\right):=ker \delta_b^j/im \delta_b^{j-1}$ with the 
pseudo-norm $\parallel .\parallel$ induced by $\parallel .\parallel_\infty$.

When dealing with a pair of multicomplexes $\left(K,L\right)$, we will 
have to
distinguish between $H_b^*\left(K,L\right)$ and 
$H_b^*\left(\mid K\mid,\mid L\mid\right)$. 

Given a bounded singular cochain $c\in C_b^j\left(\mid K\mid,\mid L\mid\right),$
we can define an antisymmetric bounded singular cochain $alt\left(c\right)\in
C_b^j\left(\mid K\mid,\mid L\mid\right)$
by
$$alt\left(c\right)\left(\sigma\right):=\frac{1}{\left(j+1\right)!}
\sum_{\pi\in S_ {j+1}}sign\left(\pi\right) c\left(\sigma\circ f_\pi\right),$$
where $S_{j+1}$ 
is the symmetric group on $e_0,\ldots,e_j$ (the 
vertices of the standard simplex) and $\sigma\circ f_\pi$ is the concatenation of $\sigma$ with the unique affine mapping $f_\pi:\Delta^j\rightarrow\Delta^j$ satisfying $f_\pi\left(e_i\right)=e_{\pi\left(i\right)}$ for $i=0,\ldots,j$. This antisymmetric bounded singular cochain
restricts to an antisymmetric bounded {\em simplicial} cochain, which we denote by
$h\left(c\right)\in C_b^*\left(K,L\right)$.
Thus, a homomorphism $h:C_b^*\left(\mid K\mid,\mid L\mid\right)\rightarrow C_b^*\left(K,L\right)$
is defined by $$h\left(c\right):=alt\left(c\right)\mid_{C_*^{simp}\left(K,L\right)}.$$
 
\begin{pro}\label{iso1} If $K$ is a minimally complete multicomplex with infinitely
many vertices in each path-component, and if $L$ is a minimally complete submulticomplex of $K$
with infinitely many vertices in each path-component, then
$$h^*:H^*_b\left(\mid K\mid,\mid L\mid\right)\rightarrow
H^*_b\left(K,L\right)$$ is an isometric isomorphism. \end{pro}

\begin{pf} For $L=\emptyset$ this is the {\em Isometry Lemma} \cite[Page 43]{gro}.
In its proof in \cite{gro} it is shown
that $h^*:
C^*_b\left(\mid K\mid\right)\rightarrow
C^*_b\left(K\right)$ has a chain homotopy inverse $g$ of norm $\parallel g\parallel\le 1$. 

If $L$ is a minimally
complete submulticomplex of $K$ then, by the construction in \cite{gro},
$g$ maps $C_b^*\left(L\right)$ to $C_b^*\left(
\mid L\mid\right)$. 
The five lemma applied to

\begin{xy}
\xymatrix{H_b^*\left(\mid L\mid\right)\ar[r]\ar[d]^{h^*}]& H_b^*\left(\mid K\mid\right)\ar[r]\ar[d]^{h^*}& H_b^*\left(\mid K\mid,\mid L\mid\right)\ar[r]\ar[d]^{h^*}& H_b^{*-1}\left(\mid L\mid\right)\ar[r]\ar[d]^{h^*}&H_b^{*-1}\left(\mid K\mid\right)\ar[d]^{h^*}\\
H_b^*\left(L\right)\ar[r]&H_b^*\left(K\right)\ar[r]&H_b^*\left(K,L\right)\ar[r]&H_b^{*-1}\left(L\right)\ar[r]&H_b^{*-1}\left(K\right)}\\
\end{xy}
\noindent shows that $h^*:H^*_b\left(\mid K\mid,\mid L\mid\right)\rightarrow
H^*_b\left(K,L\right)$ and $g^*:H^*_b\left( K, L\right)\rightarrow
H^*_b\left(\mid K\mid,\mid L\mid\right)$ are inverse isomorphisms.

The isomorphisms $h^*$ and $g^*$ have norm $\le 1$ since they are induced by cochain maps of norm $\le 1$. The inequality
$1=\parallel id\parallel \le\parallel h^*\parallel\parallel g^*\parallel\le 1$ implies that $h^*$ is an isometry.\end{pf} \\

For an $n$-dimensional compact, connected, oriented manifold, with fundamental class $$\left[M,\partial M\right]\in H_n\left(M,\partial M;{\Bbb R}\right),$$
let $$\beta_M\in 
H^n\left(M,\partial M;{\Bbb R}\right)\cong {\Bbb R}$$
be the unique class
such that $$<\beta_M,\left[M,\partial M\right]>=1.$$
By \cite[Section 1.1]{gro} $$\parallel M,\partial M\parallel =
\frac{1}{\parallel \beta_M\parallel},$$
where $\parallel \beta_M\parallel$ is the infimum of $\parallel\beta\parallel$ over all $\beta\in H^n_b(M,\partial M)$ in the preimage of $\beta_M\in H^n(M,\partial M;{\Bbb R})$.
In particular,  $$\parallel M,\partial M\parallel=0\Longleftrightarrow im\left(H^n_b
\left(M,\partial M;{\Bbb R}\right)\rightarrow H^n\left(M,\partial M\right)\right)=0.$$
 
If a continuous mapping $f:X_1\rightarrow X_2$ induces isomorphisms $\pi_0X_1\rightarrow\pi_0X_2$ and $\pi_1\left(X_1,x\right)\rightarrow\pi_1 \left(X_2,f\left(x\right)\right)$ for each $x\in X_1$, then it induces an isometric isomorphism 
of bounded cohomology (\cite[Page 40]{gro}, see also \cite[Theorem 4.1]{iva}).
If, moreover, 
$Y_1\subset X_1, Y_2\subset X_2$, and the restriction 
$$f\mid_{Y_1}:Y_1\rightarrow Y_2$$ 
induces isomorphisms $\pi_0Y_1\rightarrow\pi_0Y_2$ and $\pi_1\left(Y_1,x\right)\rightarrow\pi_1 \left(Y_2,f\left(x\right)\right)$ for each $x\in Y_1$, then, by the Relative mapping theorem \cite[Section 4.1]{gro}, $f$ induces an isometric isomorphism
$$f^*:H_b^*\left(X_2,Y_2\right)\rightarrow H_b^*\left(X_1,Y_1\right).$$
In particular, one has the following well-known fact.                        
\begin{lem}\label{degree}
If $M$ and $N$ are compact, connected manifolds 
of the same dimension and a continuous map
$f:\left(M,\partial M\right)\rightarrow\left(N,\partial N\right)$ induces isomorphisms $
\pi_1\left(M,x\right)\rightarrow \pi_1\left(N,f\left(x\right)\right)$ and $\pi_1\left(\partial M,x\right)\rightarrow\pi_1
\left(\partial N,f\left(x\right)\right)$ for each $x\in\partial M$, 
then $\parallel M,\partial M\parallel=deg\left(f\right)\parallel N,\partial N\parallel$.\end{lem}
\subsection{Aspherical multicomplexes}
\begin{pro}\label{iso2} (\cite[Section 3.3]{gro}) Let $X$ be a topological space. Then there exists an aspherical, minimally complete multicomplex
$K\left(X\right)$ such that there is an isometric isomorphism
\footnotemark \footnotetext[2]{It will be convenient for us to use notation different from Gromov's, since we will make use
of a certain functoriality of $K\left(X\right)$. So one should be aware that our $K\left(X\right)$ corresponds
to $\widetilde{K}/ \Gamma_1$ in Gromov's notation, and that our $\hat{K}\left(X\right)$ (to be defined below)
will correspond to Gromov's $\widetilde{K}$.}
$$I^*: H_b^*\left( K\left(X\right)\right)\rightarrow H^*_b\left(X\right).$$ \end{pro}

\noindent
A proof of Proposition 2 is given in \cite[Page 45-46]{gro}. Because we will need it in the proof of \hyperref[Theorem 1]{Theorem \ref*{Theorem1}}, 
we recall the construction of $\hat{K}\left(X\right)$ and $K\left(X\right)$ from \cite{gro}.\\

{\bf Geometrical description of $\hat{K}\left(X\right)$:} \\
For a topological space $X$,
$\hat{K}\left(X\right)\subset S_*\left(X\right):=\hat{S}_*^{sing}\left(X\right)$ is the multicomplex defined as follows. \\
{\em $0$-skeleton}: $\hat{K}\left(X\right)_0:=S_0\left(X\right)=\left\{x: x\in X\right\}$. \\
{\em $1$-skeleton}: For $x,y\in X$, 
$I_{x,y}$ is the set of homotopy classes rel.\ $\left\{0,1\right\}$
of maps $f:\left[0,1\right]\rightarrow X$ with $f\left(0\right)=x,f\left(1\right)=y.$
By the axiom of choice, we can choose one element $f_c$ in each homotopy class $c\in I_{x,y}$.
Define 
$$\hat{K}\left(X\right)_1:=\left\{f_c:x\not= y\in X, c\in I_{x,y}\right\}.$$
{\em $n$-skeleton}: 
Assume by recursion the $(n-1)$-skeleton  of $\hat{K}\left(X\right)$ already be defined.
For $x_0,\ldots,x_n\in X$ let $I_{x_0,\ldots,x_n}$ be the set of homotopy classes relative $\partial\Delta^n$
of maps $f:\Delta^n\rightarrow X$ with $$f\left(v_0\right)=x_0,\ldots,f\left(v_n\right)=x_n \mbox{\ and\ } f\mid_{\partial_i
\Delta^n}\in \hat{K}\left(X\right)_{n-1}\mbox{\ for\ }i=0,\ldots,n.$$
By the axiom of choice, we can choose one element $f_c$ in each class  $c\in I_{x_0,\ldots,x_n}$.

Then define
$$\hat{K}\left(X\right)_n:=\left\{f_c: x_0,\ldots,x_n\in 
X\mbox{\ pairwise\ distinct\ }, c\in I_{x_0,\ldots,x_n}\right\}.$$
Note that, by construction, $\hat{K}\left(X\right)\subset S_*\left(X\right)$ is closed under the boundary operator: $f\in \hat{K}\left(X\right)_n$ implies $\partial_if\in\hat{K}\left(X\right)_{n-1}$ for $i=0,\ldots,n$.

We have a canonical continuous mapping $$S:\mid\hat{K}\left(X\right)\mid\rightarrow X$$ which maps each chosen simplex to its image in 
$X$. According to \cite{gro} (Corollary on Page 45), $S$ is a weak homotopy equivalence. 
This implies that $S_*:C_*^{sing}\left(\mid\hat{K}\left(X\right)\mid\right)
\rightarrow C_*^{sing}\left(X\right)$ is a chain homotopy equivalence.
We denote by $T$ its chain homotopy inverse.

$T$ can be constructed by adapting the construction in \cite[Theorem 9.5]{may}.
We remark that $T$ is not induced by a continuous mapping and is not injective, but that it 
maps simplices to simplices and thus has norm $\parallel T\parallel\le 1$. 
(To compare the notation: Our $\hat{S}_*^{sing}\left(X\right)$ takes the role of the Kan complex $K$ and our $\hat{K}\left(X\right)$ corresponds to the minimal subcomplex M in \cite[Section 9]{may}. One should notice that
\cite{may} uses a slightly different construction of the minimal subcomplex $M$, namely $M$ is assumed to contain 
only one $0$-simplex in each path-component of $K$. Of course, to dispense with this assumption only simplifies the proof of \cite[Theorem 9.5]{may}.)

In fact, $T$ just sends each singular simplex $\sigma$ to the unique simplex in $\widehat{K}\left(X\right)$ which is homotopic rel.\ vertices to $\sigma$.\\


{\bf Geometrical description of $K\left(X\right)$:} \\
The multicomplex $K\left(X\right)$ is defined as the quotient $$K\left(X\right):=\hat{K}\left(X\right)/\sim$$ 
under the equivalence relation 
$$\sigma\sim\tau\in \hat{K}\left(X\right)\Longleftrightarrow \sigma_1=\tau_1.$$
That is, simplices in $\hat{K}\left(X\right)$ are identified if and only if they have the same $1$-skeleton.

The multicomplex $K\left(X\right)$ is minimally complete and aspherical by construction. \\

{\bf Pairs of spaces:}

If $Y\subset X$ is a subspace, then there is a canonical simplicial mapping $$\hat{j}:\hat{K}\left(Y\right)\rightarrow 
\hat{K}\left(X\right),$$
defined by recursion on the dimension of simplices, as follows: \\
- on the $0$-skeleton: $\hat{j}$ is defined to be the inclusion $\hat{j}:\hat{K}\left(Y\right)_0=Y\rightarrow X=\hat{K}\left(X\right)_0$,\\
- on the $1$-skeleton: for $e\in \hat{K}\left(Y\right)_1$, $\hat{j}\left(e\right)$ is the unique $1$-simplex in $\hat{K}\left(X\right)_1$ which is homotopic in $X$ rel.\ $\left\{0,1\right\}$ to $e$,\\
- ...\\
(Assume inductively, that $\hat{j}$ is defined on the $(k-1)$-skeleton such that $\hat{j}\left(\sigma\right)$ is homotopic (rel.\ vertices) to $\sigma$ for each $(k-1)$-simplex $\sigma$.)\\
- on the $k$-skeleton: Let $\sigma\in
\hat{K}\left(Y\right)_k$.  There is a homotopy (rel.\ vertices) between 
$\partial\sigma$ and $\hat{j}\left(\partial\sigma\right)$. This homotopy extends to a homotopy between 
$\sigma$ and some simplex $\sigma^\prime$ with
$\partial\sigma^\prime=\hat{j}\left(\partial\sigma\right)$. 
Define $\hat{j}\left(\sigma\right)$
to be the unique $k$-simplex in $\hat{K}\left(X\right)_k$ which is homotopic in $X$ rel\ $\partial\Delta^k$ to $\sigma^\prime$.\\

   $\hat{j}\left(\sigma\right)$ is the unique 
$k$-simplex in $\hat{K}\left(X\right)_k$ which is homotopic in $X$ rel\ $\partial\Delta^n$ to $\sigma$.\\
- ...\\

$\hat{j}$ factors to a well-defined simplicial mapping $$j:K\left(Y\right)\rightarrow K\left(X\right).$$

If $\pi_1\left(Y,y\right)\rightarrow\pi_1\left(X,y\right)$
is injective for all $y\in Y$, then
$j:K\left(Y\right)\rightarrow K\left(X\right)$ is injective. 
Indeed, if two $1$-simplices
$\gamma_1,\gamma_2\in K\left(Y\right)_1$ satisfied $j\left(\gamma_1\right)=j\left(\gamma_2\right)$,
then the representing $1$-simplices $\gamma_1,\gamma_2:\left[0,1\right]\rightarrow Y$ were homotopic rel.\ $\left\{0,1\right\}$ in $X$.
Thus $\gamma_1*\overline{\gamma}_2\in ker\left(\pi_1Y\rightarrow\pi_1
X\right)=1$, hence $\gamma_1\sim \gamma_2$ rel.\ $\left\{0,1\right\}$ in $Y$, which implies 
$\gamma_1=\gamma_2$ by construction of $K\left(Y\right)_1=\hat{K}\left(Y\right)_1$. This shows that $j$ is injective on the $1$-skeleton.
By asphericity of $K\left(Y\right)$, simplices in $K\left(Y\right)$ are determined by their $1$-skeleton and thus $j$ is injective.\\

\begin{pro}\label{iso3}(\cite[Section 4.1]{gro}, \cite[Theorem 4.1]{iva}) Let $\left(X,Y\right)$ be a pair of spaces. \\
i) If $\pi_1Y\rightarrow \pi_1X$ is injective (for each path-component of $Y$), then there exists an isometric isomorphism
$$I^*: H^*_b\left(K\left(X\right),K\left(Y\right)\right)\rightarrow
H^*_b\left(X,Y\right).$$ 
ii) If $Y$ has a decomposition $Y=Y_0\cup
Y_1$ such that 
$\pi_1Y_1\rightarrow \pi_1X$ is injective (for each path-component of $Y_1$), 
and if we let $G=\Pi\left(K\left(Y_0\right)\right)$ act on $K\left(X\right)$ 
as defined in Section 1.5.,
then there exists homomorphism
$$J^*: H_b^*\left(X,Y\right)\rightarrow H^*_b\left(K\left(X\right),GK\left(Y_1\right)\right)$$
with $\parallel J\parallel \le 1$.
 \end{pro}  
\begin{pf} 

i) We will identify $j\left(K\left(Y\right)\right)$ with $K\left(Y\right)$.
By \hyperref[iso1]{Proposition \ref*{iso1}}, $$h^*:H_b^*\left(\mid \hat{K}\left(X\right)\mid,\mid p^{-1}\left(K\left(Y\right)
\right)\mid\right)
\rightarrow H_b^*\left(\hat{K}\left(X\right),p^{-1}\left(K\left(Y\right)\right)\right)$$ is an isometric isomorphism.

The continuous mapping $S:\mid\hat{K}\left(X\right)\mid\rightarrow X$ maps $\mid p^{-1}\left(K\left(Y\right)\right)\mid$ to $Y$. By the construction in \cite[Section 9]{may}, its chain homotopy inverse $T$ maps $Y$ to $\mid p^{-1}\left(K\left(Y\right)\right)\mid$ and the chain homotopies between $ST$ and $id$ resp.\ $TS$ and $id$ preserve $Y$ resp.\ $\mid p^{-1}\left(K\left(Y\right)\right)\mid$. 
Therefore 
$$S^*: H_b^*\left(\mid\hat{K}\left(X\right)\mid,\mid p^{-1}\left(K\left(Y\right)\right)\mid\right)\rightarrow H_b^*\left(X,Y\right)$$ 
is an isomorphism with inverse $T^*$. By construction $\parallel S\parallel \le 1, \parallel T\parallel\le 1$. Together with $\parallel S^*T^*\parallel=1$ this implies that $S^*$ and $T^*$ must be isometries. 

Let $\Gamma$ be the group of simplicial automorphisms of the pair $\left( \hat{K}\left(X\right),p^{-1}\left(K\left(Y\right)\right)\right)$ and $\Gamma_i=\left\{\gamma\in\Gamma: \gamma\mid_{\hat{K}_i\left(X\right)}=id\right\}$ for $i\in{\Bbb N}$.
By \cite[Page 46]{gro}, the canonical projection
$$p:\left( \hat{K}\left(X\right),p^{-1}\left(K\left(Y\right)\right)\right)\rightarrow \left(K\left(X\right),K\left(Y\right)\right)$$ 
is a covering map with deck group $\Gamma_1$.
According to Lemma C in \cite[Section 3.3]{gro}, $\Gamma_1/\Gamma_i$ is amenable for $i\in{\Bbb N}$.

Thus, the assumptions of \hyperref[amen]{Lemma \ref*{amen}} 
(in Section 1.4 below) are satisfied for the (well-defined) action 
of $\Gamma_1/\Gamma_i$ on $\left(\hat{K}_i\left(X\right),p^{-1}\left(K_i\left(Y\right)\right)\right)$
and therefore $p^i$ has a left inverse isometry
$$Av: H_b^*\left(K_i\left(X\right),K_i\left(Y\right)\right)\rightarrow
H_b^*\left(\hat{K}_i\left(X\right),p^{-1}\left(K_i\left(Y\right)\right)\right),$$ 
for each $i\in{\Bbb N}$. In \cite[Section 3.3, Corollary D]{gro} it is claimed that the action of each $\gamma\in\Gamma$ is chain homotopic to the identity and thus $Av$ is also right inverse to $p^*$. While the first claim about chain homotopy to the identity seems not to be true, the result of $Av$ being an isometric isomorphism is indeed correct as it follows from Ivanov's work (\cite[Theorem 4.1]{iva} for the absolute case and \cite[Proposition 5.3]{bbfipp} with $\theta=1$ for the reduction from the relative to the absolute case).

Clearly, the inclusions $\left(K_i\left(X\right),K_i\left(Y\right)\right)\rightarrow
\left(K\left(X\right),K\left(Y\right)\right)$ resp.\ $\left(\hat{K}_i\left(X\right),\hat{K}_i\left(Y\right)\right)\rightarrow 
\left(\hat{K}\left(X\right),\hat{K}\left(Y\right)
\right)$ induce isomorphisms of simplicial homology and (bounded) cohomology in degrees $\le i-1$. 
Hence
$$Av: H_b^*\left(K\left(X\right),K\left(Y\right)\right)\rightarrow
H_b^*\left(\hat{K}\left(X\right),p^{-1}\left(K\left(Y\right)\right)\right)$$
is an isometric isomorphism, inverse to $p^*$. (This also would follow directly from Ivanov's work together with \cite{bbfipp}.)

Altogether, the composition $I^*=T^*\left(h^*\right)^{-1}p^*$ is an isometric isomorphism.\\

ii) 
By construction, $S$ maps $\hat{j}\left(\widehat{K}\left(Y_1\right)\right)$ to $Y_1\subset Y$. By definition of the $G$-action in Section 1.5, this implies that $S$ maps $p^{-1}\left(GK\left(Y_1\right)\right)$ to $Y$. Hence we can define $J^*:=\left(p^*\right)^{-1}h^*S^*$. By \hyperref[amen]{Lemma \ref*{amen}}, $p^*$ is an isometry. Hence $\parallel J\parallel\le 1$. \\ 

 \end{pf}



\subsection{Amenable Groups and Averaging}
\begin{df} For a group $\Gamma$ let $B\left(\Gamma\right)$ be the space of bounded real-valued functions on $\Gamma$. $\Gamma$ is called amenable if there is a $\Gamma$-invariant linear functional
$Av:B\left(\Gamma\right)\rightarrow {\Bbb R}$
such that
$inf\left(f\right)\le Av\left(f\right)\le sup\left(f\right)$ holds for all $f\in B\left(\Gamma\right)$.\end{df}
We suppose that we are given a simplicial action of
a group $G$ on a pair of multicomplexes $\left(K,L\right)$. This induces an action of $G$ on $C_b^*\left(K,L\right)$
for $g\in 
G$ and $c\in C_b^*\left(K,L\right)$, define
$gc\in C_b^*\left(K,L\right)$
by $gc\left(\left[z\right]\right):=c\left(\left[gz\right]\right)$
for $\left[z\right]\in C_*\left(K,L\right)$.
We denote
$$C_b^*\left(K,L\right)^G:=\left\{c\in C^*_b\left(K,L\right): gc=c\mbox{ for all }g\in G\right\},
\delta_i^G:= \delta\mid_{C_b^i\left(K,L\right)^G}$$ and 
$$H_b^*\left(K,L\right)^G:=ker\delta_i^G/im\delta_{i-1}^G.$$
\begin{lem}\label{amen}(\cite[Section 3]{gro}, \cite[Theorem 2.2]{iva})
If an amenable group $\Gamma$ acts on a pair of multicomplexes $\left(K,L\right)$, and $$\pi:C_b^*\left(K,L\right)^\Gamma\rightarrow
C_b^*\left(K,L\right)$$ is the
inclusion, then there is a homomorphism $$Av^*:H_b^*\left(K,L\right)\rightarrow H_b^*\left(K,L\right)^\Gamma$$ such that $$Av^*\circ \pi^*=id \mbox{\ and\ } \parallel Av^*\parallel =1.$$
\end{lem}
\begin{pf} For $L=\emptyset$ this can be proved by literally the same argument as for singular bounded cohomology in \cite[Page 39]{gro}.
The proof extends in a straightforward way to the relative case.

\end{pf}

  \subsection{Group actions on multicomplexes}     
\subsubsection{A complicated group acting on $K\left(X\right)$}
The following definitions are taken from \cite[Page 47]{gro}.
Let $\left(X,A\right)$ be a pair of spaces and $$L:=jK\left(A\right)\subset K\left(X\right).$$
We denote by $L_1$ the $1$--skeleton of $L$, by $\Omega L$ the
set of homotopy classes rel.\ $\left\{0,1\right\}$
of continuous maps $\gamma \colon \left[0,1\right]\rightarrow \mid L\mid$ with $\gamma\left(0\right)=\gamma\left(1\right)$, and by $\Omega^*L\subset\Omega L$ the subset of nontrivial homotopy classes rel.\ $\left\{0,1\right\}$.
For $x\in A$ we denote by $c_x$ the constant loop based at $x$. 
\begin{df} $$
\Pi_X\left(A\right):=\left\{\begin{array}{c}
\gamma:A\rightarrow L_1\cup\Omega L: \\
\gamma(x)=c_x\ \mbox{for all but finitely many}\ x\in A,\\
\gamma(x)(0)=x\ \forall x, \ \gamma(x)(1)\not=\gamma(y)(1)\ \forall x\not=y,\\
\forall z\in A \ \exists x\in A\ \mbox{with}\ \gamma(x)(1)=z.
\end{array}\right\}$$
$$=\left\{\begin{array}{c}\left\{\gamma_1,\ldots,\gamma_n
\right\}
:n\in {\bf N},
\gamma_1,\ldots,\gamma_n\in L_1\cup \Omega^* L,\\
\gamma_i\left(0\right)\not=\gamma_j\left(0\right), \ \gamma_i\left(1\right)\not=\gamma_j\left(1\right)\mbox{\ for\ }i\not= j,\\
\left\{\gamma_1\left(
0\right),\ldots,\gamma_n\left(0\right)\right\}=\left\{
\gamma_1\left(1\right),\ldots,\gamma_n\left(1\right)\right\}.
\end{array}\right\}.$$\end{df}


If $\gamma,\gamma^\prime\in\Omega^*L\cup L_1$ 
we denote\footnotemark\footnotetext[3]{We follow the usual convention to define the concatenation of paths 
by $\gamma*\gamma^\prime\left(t\right)=\gamma\left(2t\right)$ if $t\le\frac{1}{2}$ 
and $\gamma*\gamma^\prime\left(t\right)=\gamma^\prime\left(2t-1\right)$ 
if $t\ge \frac{1}{2}$. Unfortunately this implies that, in order to let $\Pi\left(
L\right)$ act on $K$, we will have the multiplication in $\Pi\left(L\right)$ 
such that, for example, $\left\{\gamma\right\}\left\{\gamma^\prime\right\}=\left\{
\gamma^\prime*\gamma\right\}$. We hope that this does not lead to confusion.} 
$\gamma*\gamma^\prime$ to be the unique element of $\Omega^*L\cup L_1$ in the homotopy class 
of the concatenation. 

Then we define a multiplication on
$\Pi_X\left(A\right)$ as follows:\\
given $\left\{\gamma_1,\ldots,\gamma_m\right\}$
and $\left\{\gamma_1^\prime,\ldots,\gamma_n^\prime\right\}$,
we choose a reindexing of the unordered sets $\left\{\gamma_1,\ldots,\gamma_m\right\}$
and $\left\{\gamma_1^\prime,\ldots,\gamma_n^\prime\right\}$
such that we have: $$\gamma_j\left(1\right)
=\gamma_j^\prime\left(0\right)$$ for $1\le j\le i$ and $$\gamma_j\left(1\right)
\not=\gamma_k^\prime\left(0\right)$$ for all pairs $\left(j,k\right)$ with
$j\ge i+1, k\ge i+1$.
Moreover we permute the indices $\left\{1,\ldots,i\right\}$ such that (for some $h$ with $0\le h\le i$):\\
- for $1\le
j\le h$ we have either $\gamma_j^\prime\not=\overline{\gamma_j}\in L_1\mbox{\ or\ }\gamma_j^\prime 
\not=\gamma_j^{-1}\in\Omega^*L.$\\
- for $h<
j\le i$ we have either $\gamma_j^\prime=\overline{\gamma_j}\in L_1\mbox{\ or\ }\gamma_j^\prime                                                    =\gamma_j^{-1}\in\Omega^*L.$

With this fixed reindexing we define
$$\left\{\gamma_1,\ldots,\gamma_m\right\}
\left\{\gamma_1^\prime,\ldots,\gamma_n^\prime\right\}:=
\left\{\gamma_1^\prime*\gamma_1,\ldots,\gamma_h^\prime*\gamma_h,\gamma_{i+1},\ldots,
\gamma_m,\gamma_{i+1}^\prime,
\ldots,\gamma_n^\prime\right\}.$$
Here we use the convention that $\gamma_j*\overline{\gamma_j}$ is 'empty'. (That means in this case the product is just avoided from the set on the right hand side.)\footnotemark\footnotetext[4]{Equivalently, but perhaps less intuitively, the multiplication is defined - using the first definition of $\Pi_X(A)$ - by $\gamma\gamma^\prime(x)=\gamma^\prime(x)*\gamma(\gamma^\prime(x)(1))$ for all $x\in A$.}

This multiplication\footnotemark\footnotetext[5]{We have defined a multiplication in $\Pi\left(L\right)$: if $\gamma,\gamma^\prime\in \Pi\left(L\right)$, then $\gamma*\gamma^\prime$ belongs 
to $\Pi\left(L\right)$ because $\left\{ \gamma_1^\prime*\gamma_1\left(0\right),
\ldots,\gamma_h^\prime*\gamma_h\left(0\right),\gamma_{i+1}\left(0\right),\ldots,
\gamma_m\left(0\right),\gamma_{i+1}^\prime\left(0\right),
\ldots,\gamma_n^\prime\left(0\right)\right\}=\left\{ \gamma_1^\prime*
\gamma_1\left(1\right),\ldots,\gamma_h^\prime*\gamma_h\left(1\right), \gamma_{i+1}\left(1\right),\ldots,
\gamma_m\left(1\right),\gamma_{i+1}^\prime\left(1\right),
\ldots,\gamma_n^\prime\left(1\right)\right\}$. 
The easiest way to see this is the observation that the maps $\left(
\gamma_1^\prime\left(0\right),\ldots,\gamma_n^\prime\left(0\right)\right)\rightarrow
\left(\gamma_1^\prime\left(1\right),\ldots,\gamma_n^\prime\left(1\right)\right)$ and 
$\left(\gamma_1\left(0\right),\ldots,\gamma_m\left(0\right)\right)\rightarrow\left(\gamma_1\left(1\right),\ldots,\gamma_m\left(1\right)\right)$ are permutations of $L_0$ keeping all but finitely many vertices fixed. It is then clear that the composition of two such permutations is again a permutation.   } is independent of the chosen reindexing. 
The neutral element is given by the empty set. The
inverse to $\left\{\gamma_1,\ldots,\gamma_n\right\}$ is 
given by $\left\{\overline{\gamma_1},\ldots,\overline{\gamma_n}\right\}$,
where $\overline{\gamma}_j$ is either the unique simplex resp.\ loop in the homotopy class of $t\rightarrow\gamma_j\left(1-t\right)$.
(Indeed, the above definition yields then $h=0$, thus $\left\{\gamma_1,\ldots,\gamma_n\right\}
\left\{\gamma_1^\prime,\ldots,\gamma_n^\prime\right\}$ is the empty set.)

It is well known that concatenation of paths defines an associative operation on the set of homotopy classes of paths rel.\ boundary. (Though concatenation is not associative on the set of paths.) Therefore $*$ is an associative operation on (a subset of) $L_1\cup\Omega L$. This implies associativity of the multiplication in $\Pi_X(A)$ because (using the first definition of $\Pi_X(A)$ to keep notation simpler):
{\setlength\arraycolsep{2pt}
\begin{eqnarray*}
(\gamma\gamma^\prime)\gamma^{\prime\prime}(x)&=&\gamma^{\prime\prime}(x)*\gamma^\prime(\gamma^{\prime\prime}(x)(1))*\gamma(\gamma^\prime(\gamma^{\prime\prime}(x)(1))(1)) \\
&=&\gamma^{\prime\prime}(x)*\gamma^\prime(\gamma^{\prime\prime}(x)(1))*\gamma(\gamma^{\prime\prime}(x)*\gamma^\prime(\gamma^{\prime\prime}(x)(1))(1))=\gamma(\gamma^\prime\gamma^{\prime\prime})(x)
\end{eqnarray*}}
for all $x\in A$.
Thus we have defined a group law on $\Pi_X\left(A\right)$.\\
\\
We remark that there is an
inclusion $$\Pi_X\left(A\right)\subset map_0\left(L_0,\left[\left[0,1\right],\mid L\mid\right]_{\mid K\mid}\right),$$
where $\left[\left[0,1\right],\mid L\mid\right]_{\mid K\mid}$ is the set of homotopy 
classes (in $\mid K\mid$) rel.\ $\left\{0,1\right\}$ of maps from $\left[0,1\right]$ to $\mid L\mid$, 
and $map_0\left(L_0,\left[\left[0,1\right],\mid L\mid\right]_{\mid K\mid}\right)$ 
is the set of maps $f:L_0\rightarrow \left[\left[0,1\right],\mid L\mid\right]_{\mid K\mid}$ with \\
- $f\left(y\right)\left(0\right)=y$ for all $y\in L_0$ and \\
- $f\left(.\right)\left(1\right):L_0\rightarrow L_0$ is a bijection.\\ 
This inclusion is given by sending $\left\{\gamma_1,\ldots,\gamma_n\right\}$
to the map $f$ defined by $$f\left(\gamma_i\left(0\right)\right)=\left[\gamma_i\right]$$ for $i=1,\ldots,n$, and $f\left(y\right)=\left[c_y\right]$ 
(the constant path) for $y\not\in \left\{\gamma_1\left(0\right),\ldots,\gamma_n\left(0\right)\right\}$.

The inclusion is a homomorphism with respect to the group law defined on
$map_0\left(L_0,\left[\left[0,1\right],\mid L\mid\right]_{\mid K\mid}\right)$ 
by $\left[gf\left(y\right)\right]:=\left[f\left(y\right)\right]*\left[g\left(f\left(y\right)\left(1\right)\right)\right].$

We are going to define an action of $map_0\left(L_0,\left[\left[0,1\right],\mid L\mid\right]_K\right)$ on $K$. 
This gives, in particular, an {\bf action of $\Pi_X\left(A\right)$ on $K$}.
Let $$g\in map_0\left(L_0,\left[\left(0,1\right),\mid L\mid\right]_K\right).$$
Define $gy=g\left(y\right)\left(1\right)$ for $y\in L_0$ and $gx=x$ for $x\in K_0-L_0$. This defines the action on $K_0$.

The $1$-simplices $\sigma\in K_1$ are in one-to-one-correspondence with homotopy classes (rel.\ $\left\{0,1\right\}$)
of (nonclosed) singular $1$-simplices in $\mid K\mid$ with distinct vertices in $K_0$. 
Using this correspondence,
define, for $\sigma\in K_1$, $g\sigma$ to be the unique $1$-simplex in the homotopy class 
(rel.\ $\left\{0,1\right\}$) of $$\overline{g\left(\sigma\left(0\right)\right)}*\sigma*g\left(\sigma\left(1\right)\right).$$ 
(We have $g\sigma\in K_1$
because $\overline{g\left(\sigma\left(0\right)\right)}*\sigma*g\left(\sigma\left(
1\right)\right)$ is a singular $1$-simplex with distinct vertices.
Indeed, if both vertices of $g\sigma$ agreed, then we would have
$g\left(\sigma\left(0\right)\right)\left(1\right)=g\left(\sigma\left(1\right)\right)\left(1\right)$.
But, since $g\left(.\right)\left(1\right)$ is a bijection, this would contradict $\sigma\left(0\right)\not=\sigma\left(1\right)$.)

One checks easily that $\left(gf\right)\left(\sigma\right)=g\left(f\left(\sigma\right)\right)$ for all $g,f\in 
map_0\left(L_0,\left[\left[0,1\right],\mid L\mid\right]_{\mid K\mid}\right), \sigma\in K_1$.
Thus we defined an action 
on $K_1$. 
To define the action of $map_0\left(L_0,\left[\left[0,1\right],\mid L\mid\right]_{\mid K\mid}\right)$ on all of $K$,
we claim that for a simplex
$\sigma\in K$ with $1$-skeleton $\sigma_1$, and $g\in map_0\left(L_0,\left[\left[0,1\right],\mid L\mid\right]_{\mid K\mid}\right)$,
there exists {\em some} simplex in $K$ with 1-skeleton $g\sigma_1$. Since $K$ is aspherical, this will allow a unique extension of the group action from $K_1$ to $K$.

The claim is proved by induction on the dimension of $\sigma$. First, let $\sigma$ be a $2$-simplex and denote $g\partial_0\sigma=e_0,g\partial_1\sigma=e_1,g\partial_2\sigma=e_2$. Then the concatenation of the oriented $1$-simplices $e_2,\overline{e_1}$
and $e_0$ is homotopic to the concatenation of $\partial_2\sigma,\overline{\partial_1\sigma}$ and $\partial_0\sigma$, thus null-homotopic, therefore bounds a singular $2$-simplex. By (minimal) completeness it bounds then a (unique) $2$-simplex in $K$.

Similarly, if we already have extended the action to the $(n-1)$-skeleton and let $\sigma$ be an $n$-simplex, we denote $g\partial_i\sigma=\tau_i$ for $i=0,\ldots,n$ and we observe that the composition of the $\tau_i$ is null-homotopic, hence bounds a singular simplex and, by completeness, even a simplex in $K$ which we define to be $g\sigma$.
(The vertices of $g\sigma$ are pairwise distinct, since $g\left(.\right)\left(1\right)$ permutes vertices.)
The equality $\left(gf\right)\left(\sigma\right)=g\left(f\left(\sigma\right)\right)$ extends from the $1$-skeleton $K_1$ to $K$ because, by asphericity, $\left(gf\right)\left(\sigma\right)$ is uniquely determined by $\left(gf\right)\left(\sigma_1\right)$.
This finishes the definition of the action. 
 
\begin{lem} Let $\left(K,L\right)=(K(X),jK(A))$ be a pair 
of minimally complete, aspherical multicomplexes.\\
Let $\left\{e_1,\ldots,e_n\right\}$ and 
$\left\{e_1^\prime,\ldots,e_n^\prime\right\}$ be
n-tuples of 1-simplices in $K$, such that\\
i) there is a bijection (taking multiplicities into account) between 
the set $\left\{e_1\left(0\right),e_1\left(1\right),
\ldots,e_n\left(0\right),e_n\left(1\right)\right\}$ and the set
$\left\{e_1^\prime\left(0\right),e_1^\prime\left(1\right),
\ldots,e_n^\prime\left(0\right),e_n^\prime\left(1\right)\right\}$,\\
ii) there
are oriented 1-simplices $f_1,\ldots,f_m$ in  $L$ such that 
the vertices $f_1\left(0\right),f_1\left(1\right),\ldots,f_m\left(0\right),f_m
\left(1\right)$ of $f_1,\ldots,f_m$ are all distinct and are in bijection with the 
set of vertices
of $e_1,\ldots,e_n,e_1^\prime,\ldots,e_n^\prime$, and
for each $i=1,\ldots,n$ exactly
one of the following four possibilities holds:\\
- either $\partial_0e_i=\partial_0e_i^\prime,\partial_1e_i=\partial_1e_i^\prime,
$and $e_i^\prime \overline{e_i}$ is a contractible closed loop in $\mid K\mid$,\\
- or $\partial_0e_i=\partial_0e_i^\prime, \partial_1e_i\not=\partial_1
e_i^\prime$ (resp.\ $\partial_0e_i\not=\partial_0e_i^\prime, \partial_1e=\partial_1e^\prime$)
and there is a (unique) $f\in\left\{f_1,\ldots,f_m\right\}$
such that $\partial_0f=\partial_1e_i^\prime,\partial_1f=\partial_1e_i$ (resp.\ $\partial_0f=\partial_0e_i^\prime,\partial_1f=\partial_0e_i$)
and
$e_i^\prime \overline{e_i}f$ (resp.\ $e_i^\prime f^{-1}\overline{e_i}$)
is a contractible closed loop in $\mid K\mid$,\\
- or $\partial_0e_i\not=\partial_0e_i^\prime,\partial_1e_i
\not=\partial_1e_i^\prime$ and there are (unique) $f,f^\prime\in\left\{f_1,\ldots,
f_m\right\}$
such that $\partial_0f=\partial_1e_i^\prime,\partial_1f=\partial_1e_i,\partial_0f^
\prime=\partial_0e_i^\prime,\partial_1f^\prime=\partial_0e_i
$ and $e_i^\prime
{f^{\prime}}^{-1}\overline{e_i}f$ is a contractible closed loop in $\mid K\mid$.

Then there exists some $g\in\Pi_X\left(A\right)$ with $ge_j=
e_j^\prime$ for $j=1,\ldots,n$.\end{lem}

\begin{pf} One may choose $g=\left\{f_1,\ldots,f_m,\overline{f}_1,\ldots,\overline{f}_m\right\}$.\end{pf} 
\subsubsection{A crucial observation}

\begin{obs}\label{Obs1} Let $\left(K,L\right)=(K(X),jK(A))$ be a pair of minimally complete, aspherical multicomplexes.
If $\sigma$ is any simplex in $K$
with at least one edge in $L$, then there exists $g\in\Pi\left(L\right)$ with $$ g\sigma=\overline{\sigma}.$$\end{obs}
Indeed, let $e$ be the $1$-subsimplex $\sigma$ 
which is contained in $L$. Consider $$g=\left\{e,\overline{e}\right\}\in\Pi_X\left(A\right).$$
We have $ge=\overline{e}$ because $g e$ is the unique $1$-simplex in the homotopy class of $\overline{e} *
e * \overline{e}$, which is $\overline{e}$.

Since all vertices of $\sigma$ are distinct, for edges $f\not=e$ there are
the possibilities: either $f$ has no vertex in common with $e$, or $\partial_0f\in\left\{v,w\right\}$ or/and $\partial_1f\in\left\{v,w\right\}$.

Let $w=\partial_0e,v=\partial_1e$. 
If $f$ has no vertex in common with $e$, then $g f=f$. If $\partial_0f=v$ and $\partial_1\not=w$, 
then $g f$ is the unique edge with $\partial_1\left(g f\right)=
\partial_1f$ and $\partial_0\left(g f\right)=w$. If $\partial_0f=w$
and $\partial_1f\not=v$, then $g f$ is the unique edge with $\partial_1\left(g f\right)=
\partial_1f$ and $\partial_0\left(g f\right)=v$.
If $\partial_1f=v$ and $\partial_0f\not=w$, then 
$g f$ is the unique edge with $\partial_0\left(g f\right)=
\partial_0f$ and $\partial_1\left(g f\right)=w$.
If $\partial_1f=w$ and $\partial_0f\not=v$, then 
$g f$ is the unique edge with $\partial_0\left(g f\right)=
\partial_0f$ and $\partial_1\left(g f\right)=v$. 

Thus $g$ maps the $1$-skeleton of $\sigma$ to the 1-skeleton of $
\overline{\sigma}$. Hence $g\sigma=\overline{\sigma}$, since $K$ is aspherical.\\
\\
Since $\overline{\sigma}=-\sigma$ in $C_*\left(K,L\right)$, Observation 1 implies the following:
{\begin{obs}\label{obs2}  Let $\left(K,L\right)=(K(X),jK(A))$ be a pair of minimally complete, aspherical multicomplexes.
If $\sigma$ is an $n$-simplex with at least one edge in $L$,
and if
$c\in C_b^n\left(K,L\right)^{\Pi_X\left(A\right)}
$ is some {\bf 
$\Pi_X\left(A\right)$-invariant}
bounded {\bf antisymmetric} cochain, then $c\left(\sigma\right)=0.$\end{obs}

\subsection{An application of averaging}

\subsubsection{Amenability of $\Pi_X\left(A\right)$}

\begin{lem}\label{biggroup}  If $A\subset X$ is a subspace such that $im\left(\pi_1\left(A,x\right)\rightarrow\pi_1\left(X,x\right)\right)$ is amenable for all $x\in A$, then 
$\Pi_X\left(A\right)$ is
amenable.\end{lem}

\begin{pf} There is an exact sequence 
$$1\rightarrow \bigoplus_{y\in A}im\left(\pi_1\left(A,y\right)\rightarrow\pi_1\left(X,y\right)\right)\rightarrow
\Pi_X\left(A\right)\rightarrow Perm_{fin}\left(A\right)\rightarrow 1,$$
where $Perm_{fin}$ are the permutations with finite support.

It is well known that a group is amenable if any finitely generated subgroup
is amenable. All finitely supported permutations have finite order. It follows that
any finitely generated subgroup of $Perm_{fin}\left(A\right)$ is finite and
therefore amenable. Also any finitely generated subgroup of
$\bigoplus_{y\in Y}im\left(\pi_1\left(A,y\right)\rightarrow\pi_1\left(X,y\right)\right)$
is contained in a finite sum
of amenable groups and is therefore amenable. Thus
$\Pi_X\left(A\right)$ is an amenable extension of an amenable group
and, hence, is amenable.\end{pf}

\subsubsection{Modified norms }

Let $X$ be a topological space, $Y\subset X$ a subspace, and $A\subset Y$ a subspace 
which is a union of connected components of $Y$.

As usual, we have, for each $n\in{\Bbb N}$,
the quotient complex $\ldots\rightarrow C_n\left(X,Y\right)\rightarrow
C_{n-1}\left(X,Y\right)\rightarrow\ldots$ and the short exact sequences
$0\rightarrow C_n\left(Y\right)\rightarrow C_n\left(X\right)\rightarrow
C_n\left(X,Y\right)\rightarrow 0$. We denote these maps by $\partial_n:
C_n\left(X,Y\right)\rightarrow
C_{n-1}\left(X,Y\right), i_n:C_n\left(Y\right)\rightarrow C_n\left(X\right),
j_n:C_n\left(X\right)\rightarrow
C_n\left(X,Y\right)$. On the subspace of relative cycles,
we have the boundary morphism $$d_n:C_n\left(X,Y\right)
\cap ker\left(\partial_n\right)\rightarrow C_{n-1}\left(Y\right)$$ 
which is defined by $$d_n\left(\left[\sum_{i=1}^r a_i\sigma_i\right]\right)=\sum_{i=1}^r a_i \partial\sigma_i$$ 
and which induces the (well-defined) connecting morphism $H_n\left(X,Y\right)\rightarrow H_{n-1}\left(Y\right)$.\\

Let $c\in C_n\left(X,Y\right)$ with $\partial_n c=0$. Then any $b\in j_n^{-1}\left(c\right)\subset C_n\left(X\right)$
satisfies $$\partial b=d_nc\in C_{n-1}\left(Y\right).$$

We note that $C_n\left(Y\right)$ is a direct sum of $C_n\left(Y_i\right)$ 
over the connected components $Y_i$
of $Y$. In particular, 
since $A$ is a union of connected components of $Y$,
$C_n\left(A\right)$
is a direct summand of $C_n\left(Y\right)$ and there is a canonical projection $C_n\left(Y\right)\rightarrow C_n\left(A\right)$. For $c\in C_n\left(Y\right)$, we denote by $c\mid_A\in C_n\left(A\right)$ the image of this projection.

For $\epsilon\ge 0$ define a norm on $C_n\left(X,Y\right)\cap ker\left(\partial_n\right)$ by
$$\parallel c\parallel_\epsilon^A=\inf\left\{
\parallel b\parallel + \epsilon\parallel d_n c\mid_A\parallel:
b\in j_n^{-1}\left(c\right)\right\},$$
where the norms of $b$ resp.\ $d_nc\mid_A$ are the $l^1$-norms
in $C_n\left(X\right)$ resp. $C_{n-1}\left(Y\right)$.
We get an induced pseudonorm $\parallel .\parallel_\epsilon^A$ on $H_n\left(X,Y\right)$
by $$\parallel z\parallel_\epsilon^A=\inf\left\{\parallel c\parallel_\epsilon^A: c\mbox{\ represents\ }z\right\}. $$

Moreover, on $C_b^n\left(X,Y\right)$ we have the dual norm, which we will also denote 
by $\parallel.\parallel_\epsilon^A$, and which induces a pseudonorm on the bounded
cohomology $H_b^n\left(X,Y\right)$ by $$\parallel\beta\parallel_\epsilon^A:=  \inf\left\{\parallel f
\parallel_\epsilon^A: f\mbox{\ represents\ }\beta\right\}.$$

In a completely analogous manner, one defines a pseudonorm $\parallel .\parallel_\epsilon^A$
on the bounded cohomology of pairs of multicomplexes.

If $Y\subset X$ is a $\pi_1$-injective subspace, 
then the isomorphism
$$I: H^*_b\left(K\left(X\right),K\left(Y\right)\right)\rightarrow
H^*_b\left(X,Y\right)$$ 
from section 1.3 is an isometry
for the norm $\parallel.\parallel_\epsilon^A$. This is implicit in the Relative mapping theorem, \cite[Page 57]{gro}.\footnotemark\footnotetext[6]{ A different proof of this Proposition within the framework of
homological algebra (using relative
bounded group cohomology $H_b^*\left(A\rightarrow G\right)$
defined as the
cohomology of the mapping cone of the chain map,
which is induced by the inclusion $A\rightarrow G$,
between the invariants in any compatible pair of
strong relatively injective resolution of $A$ resp.\ $G$) has been given by
Hee Sok Park in \cite{park}. However she does not prove that $H_b^*\left(A\rightarrow G\right)$ is isometric to $H_b^*\left(G,A\right)$, hence it is not clear whether her results can be applied to get results about the Gromov norm.} 

\begin{pro}\label{Prop4} If $Y\subset X$ is a $\pi_1$-injective subspace and $A$ a union of connected components of $Y$ such that $\pi_1\left(A,y\right)$ is
amenable for all $y\in A$, then $$\parallel h\parallel=
\parallel h\parallel^A_\epsilon$$ for all $h\in H_*\left(X,Y\right), \epsilon\ge 0$.\end{pro}
\begin{pf} If we assume $A=Y$, then Proposition 4 is exactly the equivalence theorem in \cite[Page 57]{gro}.
To prove the claim in generality, we give a straightforward modification of Gromov's proof.

If $\beta\in H^*_b\left(X,Y\right)$ and $h\in H_*\left(X,Y\right)$
satisfy $< \beta, h > =1$, then a standard application of the Hahn-Banach Theorem shows
$\frac{1}{\parallel \beta\parallel_\epsilon^A}=\parallel h\parallel_\epsilon^A.$
Hence it suffices to show $\parallel z\parallel_\epsilon^A=\parallel z\parallel$ for $z\in H_b^*\left(X,Y\right)$.

By \hyperref[iso1]{Proposition \ref*{iso1}} and
\hyperref[iso3]{Proposition \ref*{iso3}}, 
it suffices to prove
$\parallel c\parallel^A_\epsilon
=\parallel c\parallel$ for 
$c\in C_b^*\left(K\left(X\right),K\left(Y\right)\right)$.
 
By \hyperref[amen]{Lemma \ref*{amen}}, we may assume the bounded cochain $c$ to be invariant under the action
of the amenable group $\Pi_X\left(A\right)$. Hence, we may assume
that the relative cocycle $c$ factors over $Q$, where $Q$ is the quotient
of $F_*\left(K\left(X\right)\right)/F_*\left(K\left(Y\right)\right)$ under the relations $\overline{\sigma}=-\sigma$ and $a\sigma=\sigma$ for all $a\in \Pi_X\left(A\right)$
and all simplices $\sigma$. 
We can define in an obvious way analogues of
our norms on the dual of $Q$ and we get $\parallel c\parallel=
\parallel c^Q\parallel$ and $\parallel c\parallel^A_\epsilon=
\parallel c^Q\parallel^A_\epsilon$, where $c^Q$ is $c$ considered as a
map from $Q$ to $R$.
 
But in $Q$, any simplex $\sigma$ with an edge in $A$
becomes trivial, by \hyperref[obs2]{Observation \ref*{obs2}}.

Hence, for any
relative cycle $z\in C_*\left(X,Y\right)$, the image of $\partial z\mid_A$
in $Q$ is trivial. ($\partial z\mid_A$ is the 'restriction' of $\partial z$ to $A$ which exists because
$A$ is a union of connected components of $Y$, thus $C_*\left(A\right)$ is a direct summand of $C_*\left(Y\right)$.) Therefore $\parallel c^Q\parallel$ and
$\parallel c^Q\parallel^A_\epsilon$ agree.\end{pf}
 
\begin{cor} If $M$ is a compact manifold, $A$ a union of
connected components of $\partial M$, and $im\left(\pi_1\left(A,x\right)\rightarrow\pi_1\left(M,x\right)\right)$
is amenable for all $x\in A$, then for any $\epsilon
>0$ exists a representative $z$ of $\left[M,\partial M\right]$
with $$\parallel z\parallel\le\parallel M,\partial M\parallel+\epsilon
\mbox{\ and\ } \parallel \partial z\mid_A\parallel\le\epsilon.$$\end{cor}

In recent work with Sungwoon Kim (\cite{kk}) we use a similar argument to prove\\
$\parallel M,\partial M\parallel=\parallel M\parallel$ when $\partial M$ is amenable, where $\parallel M\parallel$ denotes the Gromov norm of the fundamental class in the locally finite homology of the interior of $M$. 

\subsection{Universal coverings of multicomplexes}


Let $K$ be a minimally complete
multicomplex and $\mid \widetilde{K}\mid$ the universal covering of 
$\mid K\mid$. 
It has a uniquely defined structure as a multicomplex $\widetilde{K}$ such that the covering $\pi$ is a simplicial map.
The Homotopy Lifting Property implies that $\widetilde{K}$ is complete. 

We construct  
a mapping
$s:K\rightarrow\tilde{K}$ with
$\pi s=id$, which maps simplices to simplices.

Define $s$ arbitrarily on the 0-skeleton $K_0$. (We have then a $\pi_1K$-equivariant bijection
$K_0\times \pi_1K = \widetilde{K_0}$ which sends $\left(v,g\right)$ to $g s\left(v\right)$.)

For $v\not=w\in K_0$ there is a unique edge $e_{vw}\in\widetilde{K_1}$ with $\partial_0e_{vw}=s\left(w\right)$ and $\partial_1e_{vw}=s\left(v\right)$.

If $f$ is any edge in $K_1$ with $\partial_0f=w,\partial_1f=v$, then the concatenation
$\overline{f} * \pi\left(e_{vw}\right)$ represents an element $\gamma\in\pi_1\left(K,v\right)$.
We define $s\left(f\right)$ to be the unique edge in $\widetilde{K}_1$ with 
$\partial_0s\left(f\right)=\gamma s\left(w\right)$ and $\partial_1s\left(f\right)=s\left(v\right)$. 

We get a 
$\pi_1K$-equivariant bijection $K_1\times\pi_1K=\widetilde{K_1}$ such that
$$\partial_0\left(\sigma,g\right)=\left(\partial_0\sigma,g \gamma\right), \partial_1\left(\sigma,g\right)=\left(\partial_1\sigma,g\right),$$
with $\gamma=\pi\left(\overline{\sigma}\right)\pi\left(e_{vw}\right),$ where $\overline{\sigma}$ is $\sigma$ with the opposite
orientation and we have used
the identification $\pi_1\left(K,v\right)=\pi_1\left(K,p\right)$ coming from their identification with the deck transformation group.

Since $\widetilde{K}$ is minimally complete, there is a unique extension of $s$ from $K_1$ to $K$, commuting with the face maps except for $K_1$.

$s$ is uniquely determined by the chosen section $s:K_0\rightarrow\widetilde{K}_0$. If $M=K\cup_AL$ is a union of connected multicomplexes $K$ and $L$ along a common submulticomplex $A$, then we may fix some connected components $\widetilde{A_i}$ of
$\pi^{-1}\left(A_i\right)$ for each connected component $A_i$ of $A$
and some connected component $\widetilde{K}$ resp.\ $\widetilde{L}$
of
$\pi^{-1}\left(K\right)$ resp.\ $\pi^{-1}\left(L\right)$, such that $\widetilde{A_i}\subset\widetilde{K}\cap\widetilde{L}$.
Then we can
define $s$ such that all $s\left(v\right)$ with $v\in A_i$ (resp.\ 
$v\in K$ resp.\ $v\in L$) 
belong to the chosen component. (This can be achieved by first defining $s$ on $A$ and then extending to $K\cup L$.)

\section{Retraction in aspherical treelike complexes}
If a group $G$ acts simplicially on a multicomplex $M$, then $C_*\left(M\right)\otimes_{{\Bbb Z}G}{\Bbb Z}$ 
are abelian groups with well-defined boundary operator $\partial_*\otimes 1$, even though $M/G$ may not be a multicomplex. (An instructive example for the latter phenomenon is the action 
of $G=\Pi_X\left(X\right)$ on $K\left(X\right)$, for a topological space $X$.) In fact, $C_*\left(M\right)\otimes_{{\Bbb Z}G}{\Bbb Z}\simeq C_*\left(M\right)\otimes_{{\Bbb R}G}{\Bbb R}$ can be considered as the quotient chain complex for the $G$-action.

We recall that in $C_*\left(M\right)
\otimes_{{\Bbb Z}G}{\Bbb Z}$ the equality $\sigma\otimes 1=g\sigma\otimes 1$ holds for 
each simplex $\sigma$. In particular, if there exists some $g\in G$ 
with $g\sigma=\overline{\sigma}$, then $\sigma\otimes 1=0\in C_*\left(
M\right)
\otimes_{{\Bbb Z}G}{\Bbb Z}$, since $\overline{\sigma}=-\sigma
\in C_*\left(M\right)$. The same remark applies to relative chain complexes.
In view of \hyperref[Obs1]{Observation \ref*{Obs1}} this implies:

If $\left(M,M^\prime\right)$ is a pair of minimally complete multicomplexes. $A$ a submulticomplex and $G=\Pi\left(A\right)$, and {\bf if $\sigma\in M$ is a simplex with at least one edge in $A$}, then 
$$\sigma\otimes 1=0\in C_*\left(M\right)\otimes_{{\Bbb Z}G}{\Bbb Z}.$$
Hence, in the Proof of Lemma 5 and Lemma 6, it is sufficient to {\bf consider only simplices $\sigma$ with 
no 
edge in $A$}.

\subsection{The 'amalgamated' case}
\begin{lem}\label{amalgamated} 
Let $\left(M,M^\prime\right)$ be a pair of minimally complete, aspherical
multicomplexes.
Let $\left(K,K^\prime\right),\left(L,L^\prime\right)$ be
pairs of 
minimally complete, aspherical submulticomplexes of $\left(M,
M^\prime\right)$. 
Assume that\\
(i)
$K^\prime=K\cap M^\prime, L^\prime=L\cap M^\prime$,\\
(ii) $M_0=K_0\cup L_0, M_0^\prime=K_0^\prime\cup L_0^\prime$,\\
(iii) the homomorphisms
$\pi_1K\rightarrow\pi_1M,\pi_1L\rightarrow
\pi_1M, \pi_1K^\prime\rightarrow\pi_1M^\prime, \pi_1L^\prime\rightarrow\pi_1M^\prime
$, induced by the respective inclusions,
are injective (for each connected component),\\
(iv) $A:=K\cap L$
is 
a path-connected
 submulticomplex of $M$, \\
the homomorphisms $\pi_1A\rightarrow\pi_1K, \pi_1A\rightarrow\pi_1L$ induced by the respective inclusions,
are injective,\\
(v) the inclusion $\left(K\cup L,K^\prime\cup L^\prime\right)\rightarrow \left(M,M^\prime\right)$ 
induces isomorphisms 
$\pi_1\left( K\cup L\right)\rightarrow\pi_1\left( M\right)$
and $\pi_1\left(K^\prime\cup L^\prime\right)\rightarrow\pi_1\left(M^\prime\right)$.

Consider the simplicial
action of $G=\Pi\left(A\right)$ on $\left(M,M^\prime\right)$ defined in section 1.5.1.

Then there is a chain homomorphism $$r:C_*\left(M\right)\otimes_{{\Bbb Z}G}{\Bbb Z}\rightarrow
C_*\left(K\right)\otimes_{{\Bbb Z}G}{\Bbb Z}\oplus C_*
\left(L\right)\otimes_{{\Bbb Z}G}{\Bbb Z}$$ in degrees $*\ge2$, mapping \footnotemark\footnotetext[7]{We do not assume $A\subset M^\prime$. Therefore we have to consider the chain complex $GC_*\left(M^\prime\right)$ to define the tensor product over ${\Bbb Z}G$. Simlarly for $C_*\left(K^\prime\right)$ and $C_*\left(L^\prime\right)$.} $GC_*\left(M^\prime\right)\otimes_{{\Bbb Z}G}{\Bbb Z}$
to $GC_*\left(K^\prime\right)\otimes_{{\Bbb Z}G}{\Bbb Z}\oplus GC_*\left(L^\prime\right)\otimes_{{\Bbb Z}G}{\Bbb Z}$, such that:\\
- if $\sigma$ is a simplex in $M$, not contained in $A$, then either $r\left(\sigma\otimes 1
\right)=\tau\otimes 1\oplus 0$ or $0\oplus\tau\otimes 1$ where $\tau$ is either a simplex in $K$ or a simplex in $L$, or $\tau=0$,\\
- if $i_K:K\rightarrow M, i_L:L\rightarrow M$ are
the inclusions, and $r_K,r_L$ the compositions of $r$ with the projection from 
$C_*\left(K\right)\otimes_{{\Bbb Z}G}{\Bbb Z}\oplus C_*\left(L\right)\otimes_{{\Bbb Z}G}{\Bbb Z}$ to the
first resp.\ second summand, then $$r_{K*}\left(i_{K*}\otimes 1\right)=id,
r_{L*}\left(
i_{L*}\otimes 1\right)=id.$$\end{lem}
 
\begin{pf} 
We may w.l.o.g.\ assume that $M$ is connected. By assumptions (iv) and (ii), this implies that $K$ and $L$ are connected.
Basically, the proof will use the fact that, unlike the complete multicomplex
$\widetilde{M}$,
in the subcomplex $\widetilde{K\cup L}$ not every pair of
vertices is connected by an edge.
(Indeed, each edge in $\widetilde{K\cup L}$
projects to an edge either in $K$ or in $L$.)
This allows to define (not unique) minimizing paths between vertices, and central simplices  for tuples of vertices, which will be used for the construction of $r$.

The plan of the proof is as follows: let $\sigma$ be a simplex in $M$, 
let $\tilde{\sigma}$ be a lift to the universal cover $\tilde{M}$, and let
$v_0,\ldots,v_n\in\tilde{M}_0$ be the vertices of $\tilde{\sigma}$. To 
each pair $\left\{v_i,v_j\right\}$ we associate a family of 'minimizing' paths $\left\{p\left(\left\{a_l^{ij}\right\};\left\{h_k^{ij}\right\}\right)\right\}$ parametrized 
by vertices $a_0^{ij},\ldots,a_{m_{ij}}^{ij}\in A_0$ and by elements $h_0^{ij},\ldots,h_{n_{ij}}^{ij}$ of $\pi_1K$ or $\pi_1L$ satisfying 
conditions described below. Associated to $\left\{v_0,\ldots,v_n\right\}$ and this 
family of 'minimizing' paths,
we construct a family of 'central' simplices \\
$\left\{\tilde{\tau}\left(\left\{a_l^{ij}\right\};\left\{h_k^{ij}\right\}\right)\subset\tilde{M}\right\}$ and 
their projections $\left\{\tau\left(\left\{a_l^{ij}\right\};\left\{h_k^{ij}\right\}\right)\subset M\right\}$, which actually lie in $K$ or $L$. 
We define $r\left(\sigma\otimes 1\right)=\tau\otimes 1\in \kk\oplus\lt$.
We show then that, for all $\tau\left(\left\{a_l^{ij}\right\};
\left\{h_k^{ij}\right\}\right)$ associated to a fixed $\tilde{\sigma}$, 
and also for all simplices $\tau$
associated to
either $g\tilde{\sigma}$ with $g\in\pi_1M$ or to $\widetilde{g\sigma}$ with $g\in G$, $\tau\otimes 1$ is the same element in $\kk\oplus\lt$.
Finally we prove that $r$ is a chain map.
Thus we can define $r\left(\sigma\otimes 1\right)=
\tau\otimes 1\oplus 0$ or $0\oplus\tau\otimes 1$. (We will henceforth 
write $\tau\otimes 1$ whenever we mean one of
$\tau\otimes 1\oplus 0$ or $0\oplus\tau\otimes 1$.)

We are going construct $r:C_*\left(M\right)\otimes_{{\Bbb Z}G}{\Bbb Z}
\rightarrow
C_*\left(K\right)\otimes_{{\Bbb Z}G}{\Bbb Z}\oplus C_*\left(L
\right)\otimes_{{\Bbb Z}G}{\Bbb Z}$ and will show in the final step that $r$ induces a relative map as in the statement of Lemma 5.\\

Let $\pi:\widetilde{M}\rightarrow M$ be the universal covering map. Let
$\widehat{K}:=\pi^{-1}\left(K\right), \widehat{L}:=\pi^{-1}\left(L\right)$ 
and $\widehat{A}:=\pi^{-1}\left(A\right)$. 
Assumption (v) implies that the universal covering
$\widetilde{K\cup L}$ is a (connected)
submulticomplex
of the universal covering $\widetilde{M}$. Assumption (ii) gives that the $0$-skeleton of $\widetilde{K\cup L}$ is the whole $0$-skeleton of $\widetilde{M}$.
From assumption (iv) we have that the connected components of $\widehat{K},\widehat{L}$ are universal coverings of $K,L$. We fix components 
$\widetilde{A},\widetilde{K},\widetilde{L}$ of $\widehat{A},\widehat{K},\widehat{L}$ such that $\widetilde{A}
\subset\widetilde{K}\cap\widetilde{L}$. Assumption (v) implies that 
$\widehat{K}\cup \widehat{L}\rightarrow K\cup L$ is the universal covering, i.e.\ $\widetilde{K\cup L}=\widehat{K}\cup\widehat{L}$.\\
In Section 1.7.\ 
we choose a section $s$ of $\pi$ 
$$s:\left(K\cup L\right)_1\rightarrow\left(\widetilde{K\cup L}\right)_1=\widehat{K}_1\cup\widehat{L}_1$$
$$\sigma\rightarrow s\left(\sigma\right)=:\tilde{\sigma}.$$
For any $f\in K_1$ (resp.\ $f\in L_1$ or $f\in A_1$)
we have $\partial_1\tilde{f}=\widetilde{\partial_1f}$
and $\partial_0\tilde{f}=\gamma \widetilde{\partial_0f}$ 
with $\gamma\in\pi_1K$ (resp.\ $\gamma\in \pi_1L$ or $\gamma\in\pi_1A$). 
As a consequence, we get the following observation.\\

(A): {\em if $g\in\pi_1M=\pi_1\left(K\cup L\right)$ and $g\tilde{e}
\in\widehat{K}_1\cup \widehat{L}_1$ is a 1-simplex
with boundary points
$\partial_1\left(g\tilde{e}\right)
=h_0\widetilde{\partial_1e}$ and $\partial_0\left(g\tilde{e}\right)=
h_1\widetilde{\partial_0e}$, 
then
$g=h_0$ and 
$$h_0^{-1}h_1\in \pi_1K\mbox{\ or\ }h_0^{-1}h_1\in \pi_1L.$$}
Indeed, we have just seen this for $g=1$. The general case follows after applying $g^{-1}$ to $g\tilde{e}$.

Moreover, if $g_1\tilde{e_1}\in\widehat{K}_1\cup \widehat{L}_1$ is a 1-simplex with 
vertices $h_{01}\tilde{v}_{11}$ and $h_{11}\tilde{v}_{01}$,
and $g_2\tilde{e_2}\in\widehat{K}_1\cup \widehat{L}_1$ 
is a 1-simplex with 
vertices $h_{02}\tilde{v}_{12}$ and $h_{12}\tilde{v}_{02}$, then:\\

(B): {\em if $g_1\tilde{e_1}$ and $g_2\tilde{e_2}$ have a common vertex $h\tilde{v}=h_{11}\tilde{v}_{01}=h_{02}\tilde{v}_{12}$, 

then either
$v\in A$ or
$h_{01}^{-1}h_{11}, h_{02}^{-1}h_{12}\in \pi_1K\mbox{\ or\ } h_{01}^{-1}h_{11}, h_{02}^{-1}h_{12}\in \pi_1L.$ }\\
Indeed, if $v\not\in A$, then $v$ is not adjacent to both, edges of $K$ and edges of $L$. \\

{\bf Minimizing paths:}
Let $v_1,v_2$ 
be vertices of $\widetilde{M}$. By a {\bf path} from $v_1$ to $v_2$ 
we mean a sequence of
1-simplices $e_1,\ldots,e_r$ in $\widetilde{K\cup L}=\widehat{K}\cup\widehat{L}$ such that $v_1$ is a vertex of $e_1$, $e_i$ and $e_{i+1}$ have a vertex
in common for $1\le i\le r-1$ and, $v_2$ is a vertex of $e_r$. (It is important that we do not allow edges in $\widetilde{M}$, where actually each pair of vertices could be connected by one edge, but only edges in $\widehat{K}\cup\widehat{L}
$.)
 
Given two vertices $v_1,v_2\in \widetilde{M}_0$, they belong 
to $\widetilde{K\cup L}_0$ because of (iii) and (v), and we have represented
them as $v_1=g_1\tilde{w_1}, 
v_2=g_2\tilde{w_2}$ with $g_i\in\pi_1\left(K\cup L\right)$ 
and $w_i:=\pi\left(v_i\right)\in \left(K\cup L\right)_0$ for $i=1,2$.

Condition (iv) implies that $\pi_1\left(K\cup L\right)=\pi_1K*_{\pi_1A}\pi_1L$ is an
amalgamated product, hence $g_1^{-1}g_2$
either belongs to $\pi_1A$ or it can be decomposed as
$$g_1^{-1}g_2=h_1\ldots h_m$$ with $h_i\in 
\pi_1K-\pi_1A$ or
$h_i\in \pi_1L-\pi_1A$ and $h_i\in\pi_1K\Longleftrightarrow h_{i+1}\in\pi_1L$ for $i=1,\ldots,m-1$.
 
Such an expression is called a {\bf normal form} of $g_1^{-1}g_2$. 
If $h_1\ldots h_m$ and $h_1^\prime\ldots h_l^\prime$ 
are two normal forms of $g_1^{-1}g_2$, then necessarily $l=m$ and for $i=1,\ldots,m$ the elements
$h_i$ and $h_i^\prime$ belong to the same equivalence class modulo $\pi_1A$. 

We emphasize that normal forms are not unique, not even up to multiplication by elements of $\pi_1A$. (That is because $im\left(\pi_1A\right)$ need not be normal.)
This will make it necessary to consider families of minimizing paths. 
 
Consider vertices $g_1\tilde{w}_1$ and $g_2\tilde{w_2}$ with $w_1,w_2\in M, g_1,g_2\in\pi_1M=\pi_1\left(K\cup L
\right)$. We define minimizing pathes in $\widetilde{K}\cup\widetilde{L}$, depending on a normal form 
$g_1^{-1}g_2=h_1\ldots h_m$ and on a choice of vertices $a_i$ in $A$. The precise definition would need to distinguish 
16 cases, according to whether $w_1,w_2$ belong to $K$ and/or 
$L$, and whether $h_1$ and/or $h_m$ belong to $\pi_1K$ and/or $\pi_1L$.

Case 1: Assume that
$w_1\in K$ and $w_2\in K$, and that 
for some (hence each) normal form $g_1^{-1}g_2=h_1\ldots h_m$ we have $h_1\in\pi_1K$ and $h_m\in\pi_1K$.
A path $e_1,\ldots,e_m$ from $g_1\tilde{w}_1$ to $g_2\tilde{w}_2$ is {\bf minimizing} if 
there is a {\em normal form} $g_1^{-1}g_2=h_1\ldots h_m$ and a 
{\em set of vertices} $a_1,\ldots,a_{m-1}\in A_0$, with $a_i\not=a_{i+1}$ for $i=1,\ldots,m-2$,
such that:\\ - $\partial_1e_1=g_1\tilde{w}_1,\partial_0e_1=g_1h_1\tilde{a}_1$,\\
- $\partial_1e_i=g_1h_1\ldots h_{i-1}\tilde{a}_{i-1},\partial_0e_i=g_1h_1\ldots h_i\tilde{a}_i$
for $i=2,\ldots,m-1$,\\
- $\partial_1e_m=g_1h_1\ldots h_{m-1}\tilde{a}_{m-1},\partial_0e_m=g_1h_1\ldots
h_m\tilde{w}_2=g_2\tilde{w}_2
$.

Case 2: Assume that $w_1\in K$ and $w_2\in K$, and that 
for some (hence each) normal form $g_1^{-1}g_2=h_1\ldots h_m$ we have $h_1\in\pi
_1K$ and $h_m\in\pi_1L$.
A path $e_1,\ldots,e_{m+1}$ from $g_1
\tilde{w}_1$ to $g_2\tilde{w}_2$ is
{\bf minimizing} if 
there is a {\em normal form} $g_1^{-1}g_2=h_1\ldots h_m$ and a 
{\em set of vertices} $a_1,\ldots,a_m\in A_0$, 
with $a_i\not=a_{i+1}$ for $i=1,
\ldots,m-1$,
such that:\\ - $\partial_1e_1=g_1\tilde{w}_1,\partial_0e_1=g_1h_1
\tilde{a}_1$,\\
- $\partial_1e_i=g_1h_1\ldots h_{i-1}\tilde{a}_{i-1},
\partial_0e_i=g_1h_1\ldots 
h_i\tilde{a}_i$
for $i=2,\ldots,m$,\\
- $\partial_1e_{m+1}=g_1h_1\ldots h_m
\tilde{a}_m=g_2\tilde{a}_m,\partial_0e_m=g_1h_1\ldots
h_m\tilde{w}_2=g_2\tilde{w}_2
$.

It should now be obvious what the analogous definition in the remaining 14 cases will be.

It follows from (A) and (B), that these paths are length-minimizing in the 
sense of being exactly the paths with 
a minimum number of edges, in $\widehat{K}\cup \widehat{L}$, between $
v_1=g_1\tilde{w}_1$ and $v_2=g_2\tilde{w}_2$. (This is because 
to each path in $\widehat{K}\cup \widehat{L}$ 
corresponds a decomposition of $g_1^{-1}g_2$ into elements 
in $\pi_1K$ and $\pi_1L$, and clearly normal forms are shortest decompositions.)

Since this latter characterization depends 
only on $v_1$ and $v_2$, we conclude\footnotemark \footnotetext[8]{One should note that all allowed edges $e_i$ in the definition of 'minimizing path' indeed exist in $\widehat{K}\cup \widehat{L}$: all neighboring points project to distinct points in $K$ resp.\ $L$
and can therefore
be joined by an edge in a translate of $\widetilde{K}$ resp.\ $\widetilde{L}$, by completeness and connectedness of $\widetilde{K}$ and $\widetilde{L}$.}: for different sections $s_1$ and $s_2$, there is a bijection between the corresponding sets of minimizing pathes from $v_1$ to $v_2$.

Since $\tilde{K}$ and $\tilde{L}$ are universal coverings of
{\em minimally complete} multicomplexes,
there is at most one edge
between two vertices. The same is true for each translate of $\tilde{K}$ resp.\ $\tilde{L}$, i.e.\ for each connected component of $\widehat{K}$ resp.\ $\widehat{L}$.
Therefore a path of length $m$ in $\widehat{K}\cup \widehat{L}$
becomes uniquely determined after fixing its $m+1$ vertices.
Hence, after fixing $a_0,\ldots,a_m\in A_0$ and 
a normal form $g_1^{-1}g_2=h_1\ldots h_m$, we get a unique path, to be 
denoted $p\left(a_0,\ldots,a_m;h_1,\ldots,h_m\right)$.

We note for later reference the following obvious observations: \\
(C1) {\em Subpaths of minimizing paths are minimizing.}\\
(C2) {\em If $e_1,\ldots,e_k$ is a minimizing path, $e_k$ projects to an edge in $K$, $e_{k+1}$ projects to an edge in $L-A$, and $e_k$ and $e_{k+1}$ have a common vertex, then $e_1,\ldots,e_k,e_{k+1}$ is a minimizing path (except if $k=1, e_1$ projects to $A$).}\\

{\bf Intersection with simplices:}
For this subparagraph, we fix a
set of vertices $\left\{v_i=g_i\tilde{u}_i\right\}_{0\le i\le n}\subset\widehat{K}_0\cup \widehat{L}_0$. For this set of vertices we consider the set 
$$\bigcup_{0\le i,j\le n}P\left(i,j\right)=\left\{r_{ij}^k: 0\le i<j\le n, r_{ij}^k\in P\left(i,j\right)\right\}, $$
where $P\left(i,j\right)$ is the set of 
minimizing paths
from $v_i$ to $v_j$.
 
Now and in the following we will consider simplices $\tilde{\tau}\in
\widehat{K}\cup \widehat{L}$ with the properties that 
no edge of $\tilde{\tau}$ belongs to $\widehat{A}_1$.
Let $\tilde{\tau}$ be such an n-simplex. 

We assert
that $\left\{\tilde{\tau}\cap r_{ij}^k: r_{ij}^k\in P\left(i,j\right)\right\}$
is the full 1-skeleton of a subsimplex of $\tilde{\tau}$. 

This means, we have to check the following claim:\\
\hspace*{1in}{\em if $\left[x,y\right], \left[z,w\right]\in\tilde{\tau}_1\bigcap \cup_{ij}P\left(i,j\right)$, \\
\hspace*{0.4in} then also $\left[x,z\right], \left[x,w\right], 
\left[y,z\right]$ and $\left[y,w\right]$
(if they exist\footnotemark \footnotetext[9]{The edges in question only exist if the corresponding vertices are different.})
 belong to $\tilde{\tau}_1\bigcap \cup_{ij}P\left(i,j\right)$.}

We will prove this assertion for $\left[x,z\right]$, assuming $x\not=z$.
The proof for the other vertices is completely analogous. 

By assumption of the claim, there is some minimizing path $\left\{\ldots,e_{l-1},e_l,e_{l+1},\ldots\right\}$ from $v_i$ to $v_j$
with
$\left[x,y\right]=e_l$ for some $l$.
By the definition of minimizing paths this means that $x=g_ih_1\ldots 
h_{l-1}\tilde{p}$ with $\tilde{p}=v_i$ or $\tilde{p}
\in\widehat{A}$ 
and that $y=g_ih_1\ldots h_l\tilde{q}$ 
with $\tilde{q}\in\widetilde{A}$
or $\tilde{q}=v_j$.
The same way, there is 
a minimizing path
$\left\{\ldots,e_{l^\prime-1}^\prime,e_{l^\prime}^\prime,e_{l^\prime+1}^\prime,\ldots\right\}$ from $v_k$ to $v_m$
with 
$\left[z,w\right]=e_{l^\prime}^\prime$, and we have
$z=g_k h_1^\prime\ldots 
h_{l^\prime-1}^\prime\tilde{p^\prime}$ with $\tilde{p^\prime}=v_k$ or $\tilde{p}
\in\widehat{A}$ and $w=g_k
h_1^\prime\ldots h_{l^\prime}^\prime\tilde{q^\prime}$ 
with $\tilde{q^\prime}\in\widehat{A}$
or $\tilde{q}=v_j$.

Note that all simplices in $\widehat{K}\cup \widehat{L}$ project 
to simplices in $K$ or in $L$. 

Assume, w.l.o.g., that $\tilde{\tau}$ projects to $K$. In particular, $\left[x,y\right]=e_l$ and $\left[z,w\right]=e_{l^\prime}^\prime$ project
to $K-A$.
By the discussion preceding observation (A),
this means that 
$h_{l}$ and $h_{l^\prime}^\prime$ belong both to 
$\pi_1K-\pi_1A$.
From the definition of minimizing paths, it follows 
then that $h_{l-1},h_{l+1},h_{l^\prime-1}^\prime,
h_{l^\prime+1}^\prime$ belong to $\pi_1L-\pi_1A$ (if they exist). 
Since, by assumption on $\tilde{\tau}$, $\pi\left(\left[x,z\right]\right)$ does not belong to $A$, hence has vertices $g\widetilde{\pi\left(x\right)}$ and 
$gh\widetilde{\pi\left(z\right)} $ for an element $h\in\pi_1K - \pi_1A$,
this implies that  
$\left[x,z\right]$ is part of some minimizing path,
namely the path $$\left\{\ldots,e_{l-2},e_{l-1},\left[x,z\right],e_{l^\prime-1}^\prime,e_{l^\prime-2}^\prime\ldots\right\}.$$
(If $h_{l-1}$ does not exist,  
then either $x=v_i$ and $\left\{\left[x,z\right],e_{l^\prime_1}^\prime,\ldots\right\}$ is a minimizing path, or $\partial_1e_{l-1}=v_i=g_1\tilde{w}_1,\partial_0e_{l-1}=x=g_1\tilde{a}_1$ and  $\left\{e_{l-1},\left[x,z\right],e_{l^\prime-1}^\prime,\ldots\right\}$ is a minimizing path. Similarly if $h_{l^\prime-1}^\prime$ does not exist.)\\


{\bf 'Central' simplices:}
We are given vertices $v_0=g_0\tilde{u_0},\ldots,v_n=g_n\tilde{u_n}\in\tilde{M}_0$, $n\ge 2$.\\
{\em We claim:} if we fix, for each index pair $\left(i,j\right)$ with $0\le i\not=j\le n$,\\
\hspace*{0.5in}a normal form
$g_i^{-1}g_j=h_1^{ij}\ldots h_{m_{ij}}^{ij}$,
vertices $a_0^{ij},\ldots,a_{m_{ij}-1}^{ij}$ and (if
necessary) $a_{m_{ij}}^{ij}\in A_0$, \\
\hspace*{0.5in}and the corresponding
minimizing path $r_{ij}\in P\left(i,j\right)$,\\
\hspace*{0.5in}then there is {\em at most one}\footnotemark \footnotetext[10]{In fact, for such an $n$-simplex to exist, the $h_{ij}$'s as well as the $a_{ij}$'s have to satisfy obvious compatibilty conditions. We will however not make use of these explicit conditions in our proof.} 
{\bf $n$-dimensional} simplex $\tilde{\tau}\in \widehat{K}\cup \widehat{L}$, 
with no edge in $\widehat{A}$,\\
\hspace*{0.5in}such that the intersection of $\tilde{\tau}$
with $\cup_{0\le i,j\le n}r_{ij}$ is the $1$-skeleton\\
\hspace*{0.5in}of an {\bf $n$-dimensional} simplex,
i.e., is the full $1$-skeleton of $\tilde{\tau}$.

A simplex $\tilde{\tau}$ as in the claim will be called the {\bf central simplex}
to the vertices $v_0,\ldots,v_n$, the normal forms 
$\left\{g_i^{-1}g_j=h_1^{ij}\ldots h_{m_{ij}}^{ij}: 0\le i<j\le n\right\}$ and the points $\left\{a_l^{ij}\right\}$.
 
We prove the claim about uniqueness of the central simplex (when it exists).
Assume there are                               
two such simplices $\tilde{\tau}_1\not=\tilde{\tau}_2$ with $dim\left(\tilde{\tau}_1\right)=dim\left(\tilde{\tau}_2\right)=n$, 
and that no edge of $\tilde{\tau}_1$ nor $\tilde{\tau}_2$ 
belongs to $\widehat{A}_1$. 
By asphericity of $\widetilde{K}$ resp.\ $\widetilde{L}$ (resp.\ their translates),
it suffices to show that $\tilde{\tau}_1$ and $\tilde{\tau}_2$ have the same $1$-skeleta. 

In the following, 'minimizing path' will mean the unique minimizing path with respect to our fixed choice of normal forms and of vertices in $A$. 
We will use the following fact: each edge of $\tilde{\tau}_1$ (resp.\ $\tilde{\tau}_2$) is contained in the minimizing path 
$r_{ij}$ from $v_i$ to $v_j$ for a unique pair 
$\left\{i,j\right\}$ of indices. This is true by a counting argument: there are 
$\frac{n\left(n+1\right)}{2}$ minimizing paths and $\frac{n\left(n+1\right)}{2}$ edges of $\tilde{\tau}_1$, each edge belongs to some minimizing path by assumption, and no minimizing path can have two consecutive edges projecting both to $K$ or both to $L$, by the 
definition of normal forms. (Note that no minimizing path hits a vertex twice, by the uniqueness of normal forms up to multiplication with elements of $\pi_1A$.)\\
For each $k$ and $l$, the minimizing path from $v_k$ to $v_l$ passes through $\tilde{\tau}_1$ as well as
through $\tilde{\tau}_2$. 
Let $\left[w_k^1,w_l^1\right]$ resp.\ $\left[w_k^2,w_l^2\right]$ be the 
intersections of this minimizing path with $\tilde{\tau}_1$ resp.\ $\tilde{\tau}_2$.

We claim: all minimizing paths from $v_k$ to some $v_i$ with $i\not=k$ pass through $w_k^1$ {\em and} $w_k^2$. Let $g_k^{-1}g_i=h_1h_2\ldots$ be the normal form corresponding to a minimaizing path from $v_k$ to $v_i$.
To prove the claim, note that, by (C1), the subpath from $v_k$ to $w_k^1$ is minimizing,
that is, the corresponding sequence $h_1,\ldots,h_s$ is a normal form (for $\Pi_{i=1}^sh_i$) with 
$h_s\in\pi_1L$ if $\tau_1$ projects to $K$ or vice versa. 
It follows from the definition of normal forms 
that also $\Pi_{i=1}^{s+1}h_i$ is a normal form, where $h_{s+1}\in\pi_1K$ 
is the element corresponding to an edge of $\tilde{\tau}$.
Since normal forms are unique up to multiplication of the $h_i$'s with elements of $\pi_1A$, we conclude that there is no minimizing path
from $v_k$ to some vertex of $\tau_1$ which does not pass through $w_k^1$.
(Otherwise, if such a minimizing path passed through the vertex $w^\prime\not=w_k^1$, then the edge $\left[w_k^1,w^\prime\right]$ would represent 
an element of $\pi_1A$ and, in particular, belong to $\widehat{A}_1$. 
By assumption, no edge of $\tilde{\tau}_1$ belongs to $\widehat{A}_1$.)

By the same argument,
all minimizing paths from $v_k$ to some vertex of $\tau_2$ have to pass through $w_k^2$. In particular, they have to contain the {\em unique} 
minimizing path from $w_k^1$ to $w_k^2$ as a 
subpath, by (C1). This, in turn, implies that all minimizing paths from $v_k$ to any $v_i, i\not=k$, contain the minimizing path from $w_k^1$ to $w_k^2$ and, in particular, contain the same edge
of $\tilde{\tau}_1$. But, by the above counting argument, it may not happen 
that several minimizing paths
pass through the same edge of $\tilde{\tau}_1$. 
This gives the contradiction.

Since it will be used later for proving
that $r$ is a chain map, we write down the following observation, which we have just proved.\\
(D): {\em If $v_0,\ldots,v_n\in\widetilde{M}_0$ are vertices of $M$, then their central simplex $\tilde{\tau}$ to 

some choice of $\left\{h_k^{ij}\right\},\left\{a_l^{ij}\right\}$, if it exists, has vertices $w_0,\ldots,w_n$ such that 

for 
any $i$ and $j$ the minimizing path from $v_i$ to $v_j$ passes through $w_i$
and $w_j$.}\\
\\
{\em Definition of $\tau$.} We denote by $\tau\left(\left\{a_l^{ij}\right\};\left\{h_k^{ij}\right\}\right)$ the projection of $\tilde{\tau}$ to $K\cup L$.\\

{\bf Well-definedness of central simplices:}
Next we claim: if we still are given $v_0=g_0\tilde{w_0},\ldots,v_n=g_n\tilde{w_n}\in\tilde{M}_0$, but vary 
$a_0,\ldots,a_m\in A_0$ and the normal forms $g_i^{-1}g_j=h_1^{ij}\ldots h_{m_{ij}}^{ij}$,
and determine $\tau:=\tau\left(\left\{a_l^{ij}\right\};\left\{h_k^{ij}\right\}\right)$, then all $\tau\otimes 1$ 
agree in
$\kk\oplus\lt$.

Let us consider $\tau\left(\left\{a_l^{ij}\right\};\left\{h_k^{ij}\right\}\right)$
and $\tau\left(\left\{{a_l^{ij}}^\prime\right\};\left\{h_k^{ij}\right\}\right)$. 
The same argument which showed uniqueness of $\tilde{\tau}\left(\left\{a_l^{ij}\right\};\left\{h_k^{ij}\right\}\right)$ lets us conclude that,
representing the vertices of $\tilde{\tau}\left(\left\{a_l^{ij}\right\};\left\{h_k^{ij}\right\}\right)$ as $\gamma_0\tilde{a}_0,\ldots,\gamma_n
\tilde{a}_n$ and the vertices of 
$\tilde{\tau}\left(\left\{{a_l^{ij}}^\prime\right\};\left\{h_k^{ij}\right\}\right)$ as $\gamma_0^\prime\tilde{a_0}^\prime,
\ldots,\gamma_n^\prime\tilde{a_n}^\prime$, we must have $\gamma_0=\gamma_0^\prime,\ldots,\gamma_n=\gamma_n^\prime$. Hence, 
corresponding tuples of edges of the two different $\tau$'s will satisfy the assumptions of lemma 3.

In fact we can simplify by argueing successively. So 
it suffices to consider the case that $\tilde{a}_0=\tilde{a}_0^\prime,\ldots,\tilde{a}_{i-1}=\tilde{a}_{i-1}^\prime,
\tilde{a}_i\not=\tilde{a}_i^\prime,\tilde{a}_{i+1}=\tilde{a}_{i+1}^\prime,
\ldots,\tilde{a}_n=\tilde{a}_n^\prime$ for some $i$ with $0\le i\le n$.
In this case, we conclude that $\tilde{\tau}\left(\left\{a_l^{ij}\right\};\left\{
h_k^{ij}\right\}\right)$ and 
$\tilde{\tau}\left(\left\{{a_l^{ij}}^\prime\right\};\left\{h_k^{ij}\right\}
\right)$ have $n$ vertices $w_j=\gamma_j\tilde{a}_j, j\not=i$, in common,
and that their distinct vertices $w_i=\gamma_i\tilde{a}_i$ and $w_i^\prime=
\gamma_i\tilde{a}_i^\prime$ are connected by an edge in $A$. Let $e$ be the 1-simplex in $T$ corresponding to this edge
and
let $g=\left\{\left[e\right],\left[\overline{e}\right]\right)\in G=\Pi_X\left(T\right)$. Then $g\circ \tau_1=\tau_2$, as has been observed in Lemma 3. (Note that
by asphericity this is implied if only $g\circ\tau_1$ and $\tau_2$
agree on the $1$-skeleton, which is true by Lemma 3.) 
But $g\circ \tau_1=\tau_2$ means $\tau_1\otimes
1=\tau_2\otimes 1\in\m$.

Now we consider $\tau\left(\left\{a_l\right\};\left\{h_k^{ij}\right\}\right)$
and $\tau\left(\left\{a_l\right\};\left\{{h_k^{ij}}^\prime
\right\}\right)$, where $g_i^{-1}g_j=h_1^{ij}\ldots h_{m_{ij}}^{ij}=h_1^{ij\prime}\ldots h_{m_{ij}}^{ij\prime}$ are different normal forms. Since 
we may argue successively, it suffices to consider the case that different normal forms occur for only one index pair $i,j$. 
For the same reason, it suffices to consider the case that there is 
$1\le s\le m_{ij}-1$ and $a\in\pi_1A$ such that ${h_s^{ij}}^\prime=h_s^{ij}a^{-1}, 
{h_{s+1}^{ij}}^\prime=ah_{s+1}^{ij}$ and ${h_l^{ij}}^\prime=h_l^{ij}$ otherwise. Now, 
one of the following two possibilities takes place:\\
- either
$g_ih_1^{ij}\ldots
h_s^{ij}\tilde{a}_s^{ij}$
is not a vertex
of the central simplex
$\tilde{\tau}\left(\left\{{a_l^{ij}}\right\};
\left\{h_k^{ij}\right\}\right)$; in this case we see that
$\tilde{\tau}\left(\left\{a_l^{ij}\right\};\left\{h_k^{ij\prime}\right\}\right)=\tilde{\tau}
\left(\left\{a_l^{ij}\right\};\left\{h_k^{ij}\right\}\right)$.\\
- or $g_ih_1^{ij}\ldots h_s^{ij}\tilde{a}_s^{ij}$ is
a vertex of 
$\tilde{\tau}\left(\left\{a_l^{ij}\right\};\left\{h_k^{ij}\right\}
\right)$; then we see that
$\tilde{\tau}\left(\left\{a_l^{ij}\right\};\left\{h_k^{ij\prime}\right\}\right)$ has $n$ vertices in common with 
$\tilde{\tau}\left(\left\{a_l^{ij}\right\};\left\{h_k^{ij}\right\}
\right)$, and
the remaining vertex in 
$\tilde{\tau}\left(\left\{a_l^{ij}\right\};\left\{h_k^{ij\prime}\right\}\right)$ resp.\ $
\tilde{\tau}
\left(\left\{a_l^{ij}\right\};\left\{h_k^{ij}\right\}\right)$ is 
$g_ih_1^{ij}\ldots h_s^{ij}a^{-1}\tilde{a}_s^{ij}$ resp.\ $
g_ih_1^{ij}\ldots h_s^{ij}\tilde{a}_s^{ij}$.
Thus the two distinct vertices of $\tilde{\tau}\left(\left\{a_l^{ij}\right\};\left\{h_k^{ij\prime}\right\}\right)
$ and $
\tilde{\tau}
\left(\left\{a_l^{ij}\right\};\left\{h_k^{ij}\right\}\right)$ 
are connected by an edge in $A$. Again, let $e$ be the 1-simplex 
in $T$ corresponding to this edge, then $g=\left\{\left[e\right],\left[\overline{e}\right]\right)\in G=\Pi_X\left(T
\right)$ maps 
$\tilde{\tau}\left(\left\{a_l^{ij}\right\};\left\{h_k^{ij\prime}\right\}\right)
$ to $
\tilde{\tau}
\left(\left\{a_l^{ij}\right\};\left\{h_k^{ij}\right\}\right)$.
Thus we get identical elements in $\m$.\\

{\bf The retraction map:}
We are going to define the retraction $r$. Given an $n$-simplex $\sigma\in M$, $n\ge2$, assume that $\sigma\otimes 1\not=0$. This implies by Observation 1 (and the remark at the beginning of Section 2) that $\sigma$ ha no edge in $A$. 

We choose a lift
$\tilde{\sigma}\in\tilde{M}$ which projects to $\sigma$. Its vertices
$v_0,\ldots v_n$ belong to $\widetilde{M}_0=
\widehat{K}_0\cup \widehat{L}_0$.\\
Now we run the above constructions for $v_0,\ldots,v_n$ to obtain a 
central simplex $\tilde{\tau}$ with $dim\left(\tilde{\tau}\right)=n$ and no edge in $\widehat{A}=\pi^{-1}\left(A\right)$.
 
If there is no central
simplex $\tilde{\tau}$, with $dim\left(\tilde{\tau}\right)=n$ and no edge in $\widehat{A}$, we define $r\left(\sigma\otimes 1\right)=0$.
If $\tilde{\tau}$ is a central simplex (with no edge in $\widehat{A}$) and if $\tau=\pi\left(\tilde{\tau}\right)$ is 
its projection, then its equivalence class $\tau\otimes 1$
is unique and we define $r\left(\sigma\otimes
1\right):=\tau\otimes 1\oplus 0$ if $\tau\in C_*\left(K\right)$ resp.\ 
$r\left(\sigma\otimes
1\right):=0\oplus\tau\otimes 1$ if $\tau\in C_*\left(L\right)$. (Note that $\tau$ has no edge in $\widehat{A}$, hence we can not have 
$\tau\in C_*\left(K\right)$ and $\tau\in C_*\left(L\right)$ at the same time.)
Clearly $r\left(\sigma\otimes 1\right)=\sigma\otimes 1$ if $\sigma\in C_*\left(K\right)$ or $\sigma\in C_*\left(L\right)$. 
 
We are going to show that 
$r$ is well-defined, i.e.\\
- independent of the lift $\tilde{\sigma}$ for the given $\sigma$, and\\
- $\sigma\otimes 1=\sigma^\prime\otimes 1$ implies $r\left(\sigma\otimes 1\right)=r\left(\sigma^\prime\otimes 1\right)$. 

If $\tilde{f}:\tilde{M}\rightarrow \tilde{M}$ is a simplicial self-map of the universal covering, which is a lift for a map $f:M\rightarrow M$ that
maps $K$ to $K$ and
$L$ to $L$, then $\tilde{f}$ 
maps minimizing paths from $v_i$ to $v_j$ to minimizing paths from 
$\tilde{f}\left(v_i\right)$ to $\tilde{f}\left(v_j\right)$. Hence, simplices $\tilde{\tau}$ intersecting the
family of minimizing paths associated to $v_0,\ldots,v_n$ in the full $1$-skeleton of an $n$-simplex are mapped by $\tilde{f}$ to simplices $\tilde{f}\left(\tilde{\tau}\right)$ intersecting the
family of minimizing paths 
associated to $\tilde{f}\left(v_0\right),\ldots,\tilde{f}\left(v_n\right)$ in the full 1-skeleton of an $n$-simplex.
Thus, if $\tilde{\tau}$ belongs to the family of 'central' simplices associated to some simplex
$\tilde{\sigma}$, then $\tilde{f}\left(\tilde{\tau}\right)$  belongs to the family of 'central' simplices associated to
$\tilde{f}\left(\tilde{\sigma}\right)$. 

We conclude: $r\left(\sigma\otimes 1\right)$ does not depend on the choice of $\sigma$ in its $G$-orbit,
nor on the choice of $\tilde{\sigma}$, for fixed $\sigma$, 
in the orbit of the deck transformation
group.\\

{\bf Compatibility with the $\partial$-operator:}
It remains to show that $r$ is a chain map, i.e., that $\partial r\left(\sigma\otimes 1\right)=r\left(
\partial\left(\sigma\otimes 1\right)\right)$ holds for all simplices $\sigma$ in $M$.\\
Case 1: $\sigma$ is a simplex 
with $r\left(\sigma\otimes 1\right)\not=0$.

Let $v_0,\ldots,v_n$ be the vertices of a lift $\tilde{\sigma}$, and $\left\{h_k^{ij}\right\}$, $\left\{a_l^{ij}\right\}$ such that a central simplex $\tilde{\tau}$ exists.
It is obvious, if we consider the set of vertices without $v_k$ and the corresponding restricted sets of normal forms and vertices, that we get the $k$-th face of $\tilde{\tau}$ as a central simplex.
That implies $r\left(\partial_k\sigma\otimes 1\right)=\partial_kr\left(\sigma\otimes 1\right)$ for all $k$, and hence $r\left(\partial\sigma\otimes 1\right)=\partial r\left(\sigma\otimes 1\right)$.\\
Case 2: $\sigma$ is a simplex 
with $r\left(\sigma\otimes 1\right)=0$.

If $r\left(\partial_k\sigma\otimes 1\right)=0$ for all faces $\partial_k\sigma$ of $\sigma$, we conclude $r\left(\partial\sigma\otimes 1\right)=0$.\\
So assume that for some face $\partial_k\sigma$ of $\sigma$ we have
$r\left(\partial_k\sigma\otimes 1\right)=
\tau\otimes 1$ for some $(n-1)$-simplex $\tau$.
That means that, for the vertices $v_0,\ldots,v_{k-1},v_{k+1},
\ldots,v_n$ and some choice of $\left\{h_k^{ij}\right\},\left\{a_l^{ij}\right\}$ we have the central 
simplex $\tilde{\tau}$. Here, the $h_k^{ij}$ and $a_l^{ij}$ are only chosen for $i\not=k,j\not=k$. By observation (D),
$\tilde{\tau}$ has vertices $w_0,\ldots,w_{k-1},w_k,\ldots,w_n$
such that the minimizing path from $v_i$ to $w_j$ passes through $w_i$ for all 
$i\not=k\not=j$. \\
For $i\not=k$ consider the set $P_{ik}^\prime$ of the minimizing paths 
from $v_k$ to $w_i$. For all $i\not=k$ choose some minimizing path $p_{ik}$ in $P_{ik}^\prime$. Let $u_{ik}$ be the vertex such that the edge $\left[u_{ik},w_i\right]$ is the last edge of $p_{ik}$.

Assume that for all $i\not=k$, the vertex $u_{ik}$ were not in the same path-component of $\widehat{A}$ as any of $\left\{w_0,\ldots,w_{k-1},w_{k+1},\ldots,w_n\right\}$. For any $i_1,i_2\not=k$, since $p_{i_1k}$ and $p_{i_2k}$
are minimizing paths, this implies that 
$u_{i_1k}$ and $u_{i_2}k$ must belong to the same component of $\widehat{A}$. That is, all $u_{ik}$ (with $i\not=k$) belong to the same component of $\widehat{A}$. Fix some vertex $u$ in this component of $\widehat{A}$, then we get that $u,w_0,\ldots,w_{k-1},w_{k+1},\ldots,w_n$ would span a central simplex for $v_0,\ldots,v_n$, contradicting the assumption $r\left(\sigma\otimes 1\right)=0$.

Thus there exists some $i,j\not=k$ such that $u_{ik}$ is in the same component of $\widehat{A}$ as $w_j$. This implies that $u_{lk}$ is in the same component of $\widehat{A}$ as $w_j$ for each $l$ with $l\not=j,l\not=k$. (Indeed, the last edge in $p_{ik}-\left[w_j,w_i\right]$ does not belong to $\tilde{\tau}$, therefore $
p_{ik}-\left[w_j,w_i\right]\cup \left[w_j,w_l\right]$ is a minimizing path from $v_k$ to $w_l$.)

But this means that $\tilde{\tau}$ is a central simplex for $v_0,\ldots,v_{j-1},v_{j+1},
\ldots,v_n$. Hence
the central simplices to $v_0,\ldots,v_{k-1},v_{k+1},
\ldots,v_n$ are exactly the central simplices to
$v_0,\ldots,v_{j-1},v_{j+1},
\ldots,v_n$, and the corresponding simplices have the same orientation if and only if
$j-k$ is {\em odd}. Therefore, $r\left(\partial_k\sigma\otimes 1\right)$ cancels against $r\left(\partial_j\sigma\otimes1\right)$. As we find such a $j$ for any $k$ with
$r\left(\partial_k\sigma\otimes 1\right)\not=0$, we get
$r\left(\partial \sigma\otimes 1\right)=0$.

As a conclusion, $r$ maps cycles to cycles and boundaries to 
boundaries.\\
To conclude, we have constructed a map $r$ with the following properties:\\
- for each $\sigma$ we have $r\left(\sigma\otimes 1\right)=\tau\otimes 1\oplus 0$ or 
$0\oplus\tau\otimes 1$ with $\tau $ a simplex either in $K$ or in $L$ (or $\tau=0$). Indeed 
we have constructed $\tau$ 
with the help of the condition that $\tilde{\tau}$ is a
simplex in $\widetilde{K}\cup\widetilde{L}$, hence $\tau$ is a simplex 
either in $K$ or in
$L$, \\
- if $\sigma$ has an edge in $A$, then ($\sigma\otimes 1=0$ and) $r\left(\sigma\otimes 1\right)=0$,\\
- if $\sigma$ is a simplex in $K$ with $\sigma\otimes 1\not=0$, then 
$r_K\left(\sigma\otimes 1\right)=\sigma\otimes 1$,\\
- if $\sigma$ is a simplex in $L$ with $\sigma\otimes 1\not=0$, then
$r_L\left(\sigma\otimes 1\right)=\sigma\otimes 1$.

The last two
items imply $r_{K*}\left(i_{K*}\otimes 1\right)=id$ and
$r_{L*}\left(i_{L*}\otimes 1\right)=id$.\\

{\bf Relative construction:}
Finally, we prove that $r$ defines a well-defined map on the relative chain complex $C_*\left(M,M^\prime\right)\otimes_{{\Bbb Z}G}{\Bbb Z}$, i.e.\ that $r$ maps 
$C_*\left(M^\prime\right)\otimes_{{\Bbb Z}G}{\Bbb Z}$ to
$C_*\left(K^\prime\right)\otimes_{{\Bbb Z}G}{\Bbb Z}\oplus
C_*\left(L^\prime\right)\otimes_{{\Bbb Z}G}{\Bbb Z}$. That is, we need to show that the construction of $r$ guarantees
that $r$ is mapping simplices in $M^\prime$ to
simplices in $M^\prime$.


For the
universal
covering $\pi:\widetilde{M}\rightarrow M$ denote the subcomplex
$\widehat{M}^\prime:=\pi^{-1}\left(M^\prime\right)$.
Let $\sigma$ be a simplex in $M^\prime$. Since we assume $
K\cap L=A$ to be path-connected, each
path component of $M^\prime$ either intersects $A$ or
is completely contained either in $K$ or in $L$.
So there are three possibilities for $\sigma$. Either $\sigma$ is contained in $K$, or $\sigma$ is 
contained in $L$, or $\sigma$ is contained in a path component of 
$M^\prime$ which intersects $A$. In the first two cases,
if $\sigma$ is contained in $K$ or $L$, 
then $r\left(\sigma\otimes 1\right)=\sigma\otimes 1$ and we are done.

We are left with the third case: we have that $\sigma\in M_c^\prime$, where $M_c^\prime$ denotes a path 
component of $M^\prime$ that intersects $A$.
Hence $\tilde{\sigma}$ belongs to some path component 
$\widehat{M}_c^\prime$
of $\widehat{M}^\prime$ (which may intersect several path components of $\widehat{A}=\pi^{-1}\left(A\right)$, since $\pi_1A\rightarrow\pi_1M_c^\prime$ need not be surjective).

First we note that, in the construction of Section 1.7.,
we may choose all $s\left(v\right)$, for $
v\in M_c^\prime$, in the 
same connected component $\widehat{M_c}^\prime$ and thus we may choose, by completeness,
$e_{vw}\in \left(\widehat{M_c}^\prime\right)_1$ whenever $v,w\in \left(
M_c^\prime\right)_0$.
This implies, for any edge $f\in \left(M_c^\prime\right)_1$,
that $\partial_0\tilde{f}=\widetilde{\partial_0f}$
and $\partial_1\tilde{f}=\gamma\widetilde{\partial_1f}$ for some
$\gamma\in \pi_1M_c^\prime$.

Let $v_i,v_j$ be two vertices of $\tilde{\sigma}$.
Since the edge $\left[v_i,v_j\right]$ belongs to $\widehat{M}_c^\prime$,
we have $v_i=g_i\widetilde{\pi\left(v_i\right)}, v_j=g_j\widetilde{\pi\left(
v_j\right)}$ with $g_i^{-1}g_j\in\pi_1M_c^\prime$.
In particular, we may choose a normal form $g_i^{-1}g_j=h_1\ldots h_m$
with $h_i$ in $im\left(\pi_1K^\prime\right)$ or $im\left(\pi_1L^\prime\right)$, by the second part of assumption (v).
It follows then from the
uniqueness of normal forms up to multiplication with elements of $im\left(\pi_1A\right)$ 
that all minimizing paths consist of edges in the $G$-orbits of edges in $\widehat{M}_c^\prime$. (If we change the normal form by multiplying two consecutive factors by an element 
of $h\in \pi_1A$ resp.\ its inverses $h^{-1}$, then the corresponding 
edges are acted upon by $h$ resp.\ $h^{-1}$ considered as 
elements of $G=\Pi\left(A\right)$.)

As a consequence, for $v_0,\ldots,v_n\in \widehat{M}^\prime_c$, any central simplex, if it exists, 
must be in the $G$-orbit of $\widehat{M}_c^\prime\subset
\widehat{M}^\prime$. Hence, either
$\tau\subset GM^\prime\cap K=GK^\prime$ or $\tau\subset GM^\prime\cap L=GL^\prime$.
This proves that $r$ induces a chain map $r:C_*\left(M,GM^\prime\right)\otimes_{{\Bbb Z}G}{\Bbb Z}
\rightarrow C_*\left(K,GK^\prime\right)\otimes_{{\Bbb Z}G}{\Bbb Z}
\oplus C_*\left(L,GL^\prime\right)\otimes_{{\Bbb Z}G}{\Bbb Z}$, finishing the proof of Lemma 5. \end{pf}   
\begin{cor}\label{coramal} Let $\left(M,M^\prime\right),\left(K,K^\prime\right),\left(L,L^\prime\right)$ 
satisfy all assumptions of \hyperref[amalgamated]{Lemma \ref*{amalgamated}}. Let
$\gamma_1\in H_b^p\left(K,K^{\prime}\right)^G, \gamma_2\in
H_b^p\left(L,L^{\prime}\right)^G$ be
bounded $G$-invariant
cohomology classes, $p\ge 2$. 
Then there
exists a class
$\gamma\in H^p_b\left(M,M^{\prime}
\right)^G$ satisfying
$\parallel \gamma\parallel\le max\left\{\parallel\gamma_1\parallel,
\parallel\gamma_2\parallel\right\}$, such that the restriction of $\gamma$ to
$K$ resp.\ $L$ gives $\gamma_1$ resp.\ $\gamma_2$.\end{cor}
\begin{pf} Let $c_1$ resp.\ $c_2$ be $G$-invariant
bounded cocycles representing $\gamma_1$ resp.\ $\gamma_2$. 

Recall that $r$ is well-defined as a homomorphism $$r:C_*\left(M,GM^\prime\right)\otimes_{{\Bbb Z}G}{\Bbb Z}\rightarrow
C_p\left(K,GK^\prime\right)\otimes_{{\Bbb Z}G}{\Bbb Z}/im\left(\partial\otimes 1\right)\oplus
C_p\left(L,GL^\prime\right)\otimes_{{\Bbb Z}G}{\Bbb Z}/im\left(\partial\otimes 1\right).$$
Let $r_1$ resp.\ $r_2$ be the composition of $r$ with the
projection to the first resp.\ second summand of the direct sum.
Since $c_1$ and $c_2$ are relative cocycles, they give well-defined homomorphisms
$$c_1:C_p\left(K,GK^\prime\right)\otimes_{{\Bbb Z}G}{\Bbb Z}/im\left(\partial\otimes 1\right)\rightarrow{\Bbb R},$$
$$c_2:C_p\left(L,GL^\prime\right)\otimes_{{\Bbb Z}G}{\Bbb Z}/im\left(\partial\otimes 1\right)\rightarrow{\Bbb R}.$$
(We are using that $G$-invariant p-cocycles on a chain complex $C$ can be considered as homomorphisms $C_p\otimes_{{\Bbb Z}G}{\Bbb Z}\rightarrow {\Bbb R}$.)
Hence we may define $c:C_p\left(M,GM^\prime\right)\rightarrow{\Bbb R}$
by\\
\\
- $c\left(\sigma\right):=c_1\left(r\left(\sigma\otimes 1\right)\right)$ if $r\left(
\sigma\otimes 1\right)\in  C_p\left(K\right)\otimes_{{\Bbb Z}G}{\Bbb Z}\oplus 0$,\\
- $c\left(\sigma\right):=c_2\left(r\left(\sigma\otimes1\right)\right)$ if $r\left(\sigma\otimes1\right)
\in 0\oplus C_p\left(L\right)\otimes_{{\Bbb Z}G}{\Bbb Z}$.\\

\hyperref[obs2]{Observation \ref*{obs2}} implies in particular that $c_1\left(\sigma\otimes 1\right)=c_2\left(\sigma\otimes 1
\right)=0$ whenever $\sigma\in C_p\left(K\cap L\right)$ (in fact, $c_1,c_2$ vanish on each simplex having an edge in $A=K\cap L$),
hence there is no ambiguity in defining $c$ on 
$C_*\left(K\cap L\right)$. 

Moreover,
$c$ is $G$-invariant by construction, hence it induces a homomorphism 
$$c:C_p\left(M,GM^\prime\right)\otimes_{{\Bbb Z}G}{\Bbb Z}\rightarrow{\Bbb R}.$$

All claims of the corollary are obvious except possibly that $c$ is indeed a 
relative cocycle, i.e.,
that $\delta c=0$ and that $c$ vanishes on p-simplices $\sigma$
in $GM^\prime$. 

To prove $\delta c=0$, we consider $\tau\in C_{p+1}\left(M\right)$ and assume w.l.o.g. $r\left(\tau\otimes 1\right)\in
C_{p+1}\left(K\right)
\otimes_{{\Bbb Z}G}{\Bbb Z}$. Since $c_1$ is a cocycle, 
we get $\partial c\left(\tau\right)=
c\left(\partial\tau\right)=c_1\left(\partial r\left(\tau\otimes 1\right)\right)=\delta c_1\left(r\left(\tau\otimes 1\right)\right)=0$.

To prove that $c$ vanishes on p-simplices 
in $GM^\prime$, let $\sigma\in C_p\left(M^\prime\right)$
and
assume
w.l.o.g. that $r\left(\sigma\otimes 1\right)\in 
C_*\left(GK^\prime\right)
\otimes_{{\Bbb Z}G}{\Bbb Z}$. Then $c_1\left( r\left(\sigma
\otimes 1\right)\right)=0$ 
because $c_1\in C_b^p\left(K,K^\prime\right)$ vanishes on $K^\prime$ and thus, by $G$-invariance, on $GK^\prime$. 
Therefore $c\left(\sigma\right)=
c_1\left(r\left(
\sigma\otimes 1\right)\right)
=0$. Since $c$ is $G$-invariant, this proves that $c$ vanishes on $C_p\left(GM^\prime\right)$.\end{pf}

\subsection{The 'HNN'-case}

For a multicomplex $A$ define $A\times I$ to be the following multicomplex:\\
- the $0$-skeleton is $A_0\times\left[0,1\right]$,\\
- the $1$-skeleton is $ A_1\times\left[0,1\right]\times \left[0,1\right]$, with $\partial \left(\sigma,t,s\right)=\left(\partial_0\sigma,s\right)-\left(\partial_1\sigma,t\right)$ for all $\left(\sigma,t,s\right)\in
 A_1\times\left[0,1\right]\times \left[0,1\right]$,\\
- for $n\ge 2$, an $n$-simplex with a given $(n-1)$-skeleton in $A\times I$ exists (uniquely) if and only if the projection of the $(n-1)$-skeleton to $A$ is bounding an $n$-simplex in $A$.

If $A$ is a minimally complete, aspherical multicomplex, then $A\times I$ is a minimally complete, aspherical multicomplex, and 
$\mid A\times I\mid$ is homotopy equivalent to 
$\mid A\mid\times I$.


\begin{lem}\label{hnn} 
Let $\left(M,M^\prime\right)$ be a pair of minimally complete, aspherical
multicomplexes.
Assume that 
there are 
pairs of
minimally complete, aspherical
multicomplexes $\left(K,K^\prime\right), \left(L,L^\prime\right)$,
such that 
$K\cap L $ consists of two connected
components $B, C$, which are 
minimally complete, aspherical submulticomplexes and such that $\left(L,L^\prime\right)$ is isomorphic to $\left(B\times I,\left(B\cap M^\prime\right)\times I\right)$, for an isomorphism carrying $B\cup C
$ to $B\times\left\{0,1\right\}$. (In particular, there is an isomorphism $F:B\cong C$.)

Assume that:\\
(i) $K^\prime=M^\prime\cap K, L^\prime=M^\prime\cap L$,\\
(ii) $M_0=K_0\cup L_0, M_0^\prime=K_0^\prime\cup L_0^\prime$,\\
(iii) the homomorphisms $\pi_1K\rightarrow\pi_1M,\pi_1L\rightarrow
\pi_1M, 
\pi_1K^\prime\rightarrow\pi_1M^\prime, 
\pi_1L^\prime\rightarrow\pi_1M^\prime,$
induced by the respective inclusions,
are injective (for each connected component),\\
(iv) 
the homomorphisms $\pi_1 B\rightarrow\pi_1K, \pi_1C\rightarrow \pi_1K$, induced by the respective inclusions,
are injective,\\
(v) the inclusion $\left(K\cup L,K^\prime\cup L^\prime\right)\rightarrow \left(M,M^\prime\right)$
induces isomorphisms $\pi_1\left(K\cup L\right)\rightarrow\pi_1M$
and $\pi_1\left(K^\prime\cup L^\prime\right)\rightarrow\pi_1 M^\prime$.

Consider the simplicial action of $G=\Pi\left(B\cup C\right)$ on
$\left(M,M^\prime\right)$.

Then there is a chain homomorphism $$r:C_*\left(M,GM^\prime\right)
\otimes_{{\Bbb Z}G}{\Bbb Z}\rightarrow
C_*\left(K,GK^\prime\right)\otimes_{{\Bbb Z}G}{\Bbb Z}
$$ in degrees $
*\ge2$ such that\\
- if $\sigma$ is a simplex in $M$, then $r\left(\sigma\otimes 1
\right)=\tau\otimes 1$ and $\tau$ is simplex in $K$ or $\tau=0$,\\
- $r\left(i_{K*}\otimes 1\right)=id_{C_*\left(K,GK^\prime\right)
\otimes_{{\Bbb Z}G}
{\Bbb Z}}$.\end{lem}

\begin{pf}
The proof of Lemma 6 parallels in several aspects the proof of \hyperref[amalgamated]{Lemma \ref*{amalgamated}}. We 
will therefore only emphasize those points
which are different from the proof of Lemma 5 and else refer to the proof 
given there.

We may w.l.o.g.\ assume that $K$ (and thus $M$) are connected.

It follows from the Seifert-van Kampen Theorem
that $\pi_1\left(K\cup L\right)$ is an HNN-extension of
$\pi_1K$. We consider $\pi_1K$ as subgroup of $\pi_1\left(K\cup L\right)$ and denote by $t$ the
extending element of $$\pi_1\left(K\cup L\right)=
<\pi_1K,t\mid tat^{-1}=F_*a\ \ \forall a\in im\left(\pi_1B\right)>.$$
Here $F:B\rightarrow C$ is the simplicial isomorphism obtained as restriction
from the given isomorphism $B\times I\rightarrow L$. \\
(Even though $B,C$ are disjoint, we will consider $im\left(\pi_1B\right),im\left(\pi_1
C\right)$ as subgroups of $\pi_1\left(
K\cup L\right)$ by the usual conjugation with some fixed path from the basepoint 
to $B$ resp.\ $C$. We will write $im\left(\pi_1A\right)$ for the union of 
$im\left(\pi_1B\right)$ and $im\left(\pi_1C\right)$ as a subset of $\pi_1K$.)

Let $\pi:\widetilde{M}\rightarrow M$ be the universal covering. 
We can choose connected components $\widetilde{B},\widetilde{C},
\widetilde{L}$ of $\pi^{-1}\left(B\right),\pi^{-1}\left(C\right)$ 
resp.\ $\pi^{-1}\left(L\right)$, in such a way that $\widetilde{B}
\subset\widetilde{L}$ and $\widetilde{C}\subset\widetilde{L}$. 
Moreover we can choose a connected component $\widetilde{K}$ of $\pi^{-1}\left(K\right)$ such that $\widetilde{B}\subset\widetilde{K}$ but 
$\widetilde{C}\not\subset\widetilde{K}$. It follows from (iii) and (iv) that the restrictions of $\pi$ to $\widetilde{K},\widetilde{L}$ are indeed universal covering maps.

By (v), the universal covering $
\widetilde{K\cup L}$ is a submulticomplex of $\widetilde{M}$. By (ii), $M_0=K_0\cup L_0$. Hence
(v)
implies that the 0-skeleton of
$\widetilde{K\cup L}$ equals $\widetilde{M}_0$. \\
In Section 1.7.\ 
we constructed a specific section 
$s:\left(K\cup L\right)_1\rightarrow\left(\widetilde{K\cup L}\right)_1$ of $\pi$. (We will again
denote $\tilde{\sigma}:=s\left(\sigma\right)$ for simplices $\sigma$ in $K\cup L$.)

If $v$ is a vertex in $B$ resp.\ $C$ resp.\ $L$ then we choose $\tilde{v}$ in $\widetilde{B}$ resp.\ $\widetilde{C}$ resp.\ $\widetilde{L}$. If $v$ is a
vertex in $K-C$, then we choose $\tilde{v}\in \widetilde{K}$. (Thus $\tilde{v}
\not\in\widetilde{K}$ for $v\in C\subset K$. But $\tilde{v}\in\widetilde{C}$ implies $\tilde{v}=\gamma\tilde{\tilde{v}}$ with $\tilde{\tilde{v}}\in\widetilde{K}, \gamma=th, h\in\pi_1K$.)

If $v$ and $w$ are vertices in $K-C$ (resp.\ in $L$, resp.\ in $B$ or $C$), then
we can
and will choose $e_{vw}\in\widetilde{K}$ (resp.\ $\widetilde{L}$, resp.\ $
\widetilde{B}$ or $\widetilde{C}$) in the construction of Section 1.7.
This implies that, for any edge
$f\in K_1$, with $\partial_1f\not\in C$, $ \tilde{f}$
must belong to
$\widetilde{K}$
and satisfy $\partial_1\tilde{f}=\tilde{\partial_1f}$
and $\partial_0\tilde{f}=\gamma \tilde{\partial_0f}$
with $\gamma\in\pi_1\left(K\right)\subset\pi_1\left(K\cup L\right)$.
Analogous observations apply for edges $f\in L_1$
resp.\ $f\in B_1$ or $f\in C_1$. 

Finally, if $f\in K_1$ with $\partial_1f\in C$, then $\tilde{f}$
belongs to some component of 
$\pi^{-1}\left(K\right)$, i.e.\ some $\gamma\widetilde{K},\gamma\in
\pi_1M$. Thus $\partial_0\tilde{f}=\gamma\tilde{\partial_0f}$. But $\partial_1\tilde{f}=\widetilde{\partial_1f}\in\widetilde{C}$
implies that $\gamma=th$ for some $h\in\pi_1K$.
In conclusion we have\footnotemark \footnotetext[11]{We are now in a somewhat different
situation as in the proof of \hyperref[amalgamated]{Lemma \ref*{amalgamated}}.
Namely, if we have two points $g_1\tilde{a}_1$ and $g_2\tilde{a}_2$ with
$a_1,a_2\in B$ and $g_2g_1^{-1}\not\in\pi_1B$,
then there is no edge in $\widetilde{K\cup L}$
with boundary points $g_1\tilde{a}_1$ and $g_2\tilde{a}_2$.
Moreover, there is a
path between $g_1\tilde{a}_1$ and $g_2\tilde{a}_2$
consisting of two edges if and only if $g_1^{-1}g_2=\left(t^{\pm1}h\right)^{\pm1}$
for some $h\in\pi_1K$ (and $t$ the extending element). The same
remark applies if $a_1,a_2\in C$.}:

(A): {\em if $g\tilde{e}\in\left(\widetilde{K\cup L}\right)_1$ is a 1-simplex
with 
$\partial_1\left(g\tilde{e}\right)
=h_0\widetilde{\partial_1e}$ and $\partial_0\left(g\tilde{e}\right)=
h_1\widetilde{\partial_0e}$,  \\
then
$g=h_0$ and

$h_0^{-1}h_1$ is either an element of $\pi_1L\subset \pi_1K$ (if $e\in L_1$), 
or of the form
$\left(t^{\pm 1}h\right)^{\pm1}$ with $h\in\pi_1K$.}\\
{\bf Minimizing paths:}
Given two vertices $v_0,v_1\in\widetilde{M}_0$,
we may represent them
as $v_1=g_1\tilde{w_1}, v_2=g_2\tilde{w}_2$ with $g_i\in\pi_1 M$ and $w_i\in M_0$ for $i=1,2$. $\pi_1M$ is an HNN-extension of $\pi_1K$, hence $g_1g_2^{-1}$
has an expression $g_1g_2^{-1}=h_1\ldots h_m$
with $h_i\in\pi_1K$ or $h_i= t^{\pm1}
$ such that $h_i\in\pi_1K$ if and only if $h_{i+1}=t^{\pm1}$. ({\em This time, we allow $h_i=1\in\pi_1K$ except for $i=1$ and $i=m$.})
This expression, which we will call a {\bf normal form},
is unique upon compatible changes of
the $h_i\in\pi_1K$ by multiplication 
with elements of
$im\left(\pi_1A\right):=im\left(\pi_1 B\right)\cup im\left(\pi_1C\right)$ (e.g.\ $\left(h_1\ F_*a\right)t\left(h_2\right)=\left(h_1\right)t\left(a h_2\right)$ for $a\in\m\left(pi_1B\right), h_1,h_2\in\pi_1K$).

Minimizing paths between $v_1=g_1\tilde{w}_1$ and $v_2=g_2\tilde{w}_2$ are defined in a
similar manner as in the proof of \hyperref[amalgamated]{Lemma \ref*{amalgamated}} (with the same number of case 
distinctions). The only formal difference (which is caused by the Remark above) is that the $a_i$'s have to be chosen in $B$ and $C$ alternatingly and
that a minimizing path of length 2
between
$g_1\tilde{a}_1$ and $g_1ht\tilde{a}_3$ with $a_1,a_3\in B$ and 
$h\in\pi_1K$ is either of
of the form $\left\{e_1,e_2\right\}$ with $\partial_1e_1=g_1\tilde{a}_1, \partial_0e_1=\partial_1e_2=g_1hb\tilde{a}_2, \partial_0e_2=ght\tilde{a}_3$ with $b\in im\left(\pi_1C\right)$ and $a_2\in C$. Similarly with $B,C$ interchanged.

It follows from the remark above
that these paths are length-minimizing in the sense of being exactly the paths between $v_1$ and $v_2$ with a minimal number of edges.

Let $v_0,\ldots,v_n$ be vertices of $\widetilde{K\cup L}$.
For each $i,j\in\left\{0,\ldots,n\right\}$ define $P\left(i,j\right)$ 
as in the proof of \hyperref[amalgamated]{Lemma \ref*{amalgamated}}. Let $\tilde{\tau}$ be a simplex in $\widetilde{K\cup L}$, with no edge in $\pi^{-1}\left(B\right)\cup\pi^{-1}\left(C\right)$.
Then the same proof as in the 'amalgamated' 
case shows that, for each simplex $\tilde{\tau}$, 
the set $\left\{\tilde{\tau}\cap r_{ij}^k: r_{ij}^k\in
P\left(i,j\right)\right\}$ is the full 1-skeleton of a subsimplex of 
$\tilde{\tau}$. Exactly as in the proof of \hyperref[amalgamated]{Lemma \ref*{amalgamated}} this allows us to obtain well-defined 
'central' simplices and thus to construct the retraction $r:C_*\left(M,GM^\prime\right)\otimes_{{\Bbb Z}G}{\Bbb Z}\rightarrow
C_*\left(K,GK^\prime\right)\otimes_{{\Bbb Z}G}{\Bbb Z}\oplus C_*\left(L,GL^\prime\right)\otimes_{{\Bbb Z}G}{\Bbb Z}$ in degrees $*\ge2$. 

We observe that simplices $\tau$ in $\widetilde{L}$ of dimension $\ge 2$, whose edges belong to minimizing paths, necessarily have 
some edges which project to $B\cup C$.
Since this implies $\tau\otimes 1=0$ by \hyperref[obs2]{Observation \ref*{obs2}}, we conclude that central simplices necessarily belong to $K$. Therefore we
have $r:C_*\left(M,GM^\prime\right)\otimes_{{\Bbb Z}G}{\Bbb Z}\rightarrow
C_*\left(K,GK^\prime\right)\otimes_{{\Bbb Z}G}{\Bbb Z}$ with the wished 
properties.
\end{pf}\\
\\
In an analogous manner to \hyperref[coramal]{Corollary \ref*{coramal}}, we conclude
\begin{cor}\label{corhnn} 
 Let $\left(M,M^\prime\right),\left(K,K^\prime\right),\left(L,L^\prime\right)$
satisfy all assumptions of \hyperref[hnn]{Lemma \ref*{hnn}}. \\
Let
$\gamma_1\in H_b^p\left(K^G,K^{\prime}\right)^G$ be
a bounded $G$-invariant cohomology class, $p\ge 2$.
Then
there exists a class
$\gamma\in H^p_b\left(M^G,M^{\prime}
\right)^G$, satisfying
$\parallel \gamma\parallel\le \parallel\gamma_1\parallel$,
such that restriction of $\gamma$ to $K$ gives $\gamma_1$.
\end{cor}

\subsection{Locally finite chains}

In this paper
we are interested in compact manifolds 
with boundary. 
The generalization
of \hyperref[amalgamated]{Lemma \ref*{amalgamated}} resp.\ \hyperref[hnn]{Lemma \ref*{hnn}} to locally finite chains (and disconnected intersections) will be needed in \cite{k4} and can be used to 
generalize Theorem 1 for the simplicial volume of noncompact manifolds for cut and paste along {\em compact} submanifolds.

\begin{pro}\label{pro5} Assume that $\left(K,K^\prime\right),\left(L,L^\prime\right),\left(M,M^\prime\right)$ satisfy all assumptions of 
\hyperref[amalgamated]{Lemma \ref*{amalgamated}} or \hyperref[hnn]{Lemma \ref*{hnn}},
respectively. Assume, in addition, that $\mid A\mid=\mid K\cap L\mid$ is compact. 

Then $r$ extends to a chain homomorphism (for $*\ge 2$)
$$r:C_*^{lf}\left(M\right)\otimes_{{\Bbb Z}G}{\Bbb Z}\rightarrow
C_*^{lf}\left(K\right)\otimes_{{\Bbb Z}G}{\Bbb Z}\oplus C_*^{lf}\left(L\right)\otimes_{{\Bbb Z}G}{\Bbb Z},$$ mapping $C_*^{lf}\left(GM^\prime\right)\otimes_{{\Bbb Z}G}{\Bbb Z}$
to $C_*^{lf}\left(GK^\prime\right)\otimes_{{\Bbb Z}G}{\Bbb Z}\oplus C_*^{lf}\left(GL^\prime\right)\otimes_{{\Bbb Z}G}{\Bbb Z}$ and still satisfying
summand, then 
$r_{K*}\left(i_{K*}\otimes 1\right)=id,
r_{L*}
\left(i_{L*}\otimes 1\right)=id$.
\end{pro}

\begin{pf} 
Let $\sum_{i=1}^na_i\sigma_i\otimes 1
\in C_*^{lf}\left(M,GM^\prime\right)\otimes_{{\Bbb Z}G}{\Bbb Z}$ be a locally finite chain.\\
Since $\mid A\mid$ is compact, for any
locally finite chain $\sum_{i=1}^na_i\sigma_i\otimes 1
\in C_*^{lf}\left(M,GM^\prime\right)\otimes_{{\Bbb Z}G}{\Bbb Z}$
only finitely many $\sigma_i$ intersect $A$.
If $\sigma_i$ does not intersect $A$, then $r\left(\sigma_i\otimes 1\right)=\sigma_i\otimes 1$
by construction. Thus $r\left(\sigma_i\otimes 1\right)=\sigma_i\otimes 1$ for all but finitely many $\sigma_i$ and $\sum_{i=1}^n a_ir\left(\sigma_i\otimes 1\right)$
is locally finite.
\end{pf}

\begin{pro} 
Let all assumptions of \hyperref[amalgamated]{Lemma \ref*{amalgamated}} be satisfied except that $A$ need not be path-connected.
Then there is a chain homomorphism $$r:C_*^{lf}\left(M\right)\otimes_{{\Bbb Z}G}{\Bbb Z}\rightarrow C_*^{lf}\left(K
\right)\otimes_{{\Bbb Z}G}{\Bbb Z}
$$
in degrees $*\ge2$, mapping
$C_*^{lf}\left(GM^\prime\right)\otimes_{{\Bbb Z}G}{\Bbb Z}$ to $C_*^{lf}\left(GK^\prime\right)\otimes_{{\Bbb Z}G}{\Bbb Z}$
such that\\
- if $\sigma$ is a simplex in $M$, then $r\left(\sigma\otimes 1 \right)=\kappa\otimes 1$, where either
$\kappa$ is a simplex in $K$
 or $\kappa=0$,\\
- $r\left(i_K\otimes id\right)=id_{C_*\left(K,GK^\prime\right) \otimes_{{\Bbb Z}G}{\Bbb Z}}$.
\end{pro}

\begin{pf}
Assume $M$ is connected. Let $A_1,\ldots,A_n$ be the path-components of $A$.

If $n=1$, i.e.\ if $A$ is connected,
then $K$ and $L$ must be connected and we can apply \hyperref[amalgamated]{Lemma \ref*{amalgamated}} (and compose $r$ with projection to the first summand).

If $n>1$, then $Q_1=N - B_1$ is path-connected, hence
we can apply \hyperref[hnn]{Lemma \ref*{hnn}} to construct $r_1$.
Succesively for $i=2,\ldots,n-1$, $Q_i=Q_{i-1}-B_i$ must be connected, hence we can apply \hyperref[hnn]{Lemma \ref*{hnn}},
to construct $r_i$. 

Finally, $Q_n=Q_{n-1}-B_n$ consists of two connected components, hence we can apply \hyperref[amalgamated]{Lemma \ref*{amalgamated}} (and compose with projection to the first summand) to construct $r_n$. Then define $r=r_nr_{n-1}\ldots r_1.$\end{pf}

\section{Glueing along amenable boundaries}

\subsection{Proof of Theorem 1} 

In this section, we prove superadditivity for simplicial volume of manifolds 
with boundary with respect to glueing along amenable subsets of the boundary.
A similar result for open manifolds is the Cutting-of theorem in \cite[Section 4.2]{gro}.
One should note that, at least for manifolds with boundary, the opposite inequality need not hold. As a counterexample one may glue solid tori along disks to get a handlebody.   
\begin{lem}\label{lem8} (i): Let $M_1,M_2$ be two compact, connected n-manifolds,
$A_1,A_2$ (n-1)-dimensional submanifolds of $\partial M_1$ resp.\ $\partial M_2$,
$f:A_1\rightarrow A_2$ a homeomorphism
and $M=\glu$ the glued manifold. Assume that $\pi_1\partial M\rightarrow \pi_1M$ is injective.

If $f_*$ induces an ismorphism
$$f_*:ker\left(\pi_1A_1\rightarrow\pi_1M_1\right)\rightarrow
ker\left(\pi_1A_2\rightarrow\pi_1M_2\right),  ker\left(\pi_1\partial A_1\rightarrow\pi_1M_1\right)\rightarrow
ker\left(\pi_1\partial A_2\rightarrow\pi_1M_2\right)$$ and if
$im\left(\pi_1A_1\rightarrow\pi_1M_1\right)$ is amenable,
then $$\parallel M,\partial M\parallel\ge\sima+\simb.$$
(ii):  Let $M_1$ be a compact, connected n-manifold,
no component of which is a 1-dimensional closed intervall,
$A_1, A_2$ disjoint (n-1)-dimensional submanifolds of $\partial M_1$,
$f:A_1\rightarrow A_2$ a homeomorphism
and $M=M_1/f$ the glued manifold. Assume that $\pi_1\partial M\rightarrow \pi_1M$ is injective.

If $f_*$ induces an isomorphism
$$f_*:ker\left(\pi_1A_1\rightarrow\pi_1M_1\right)\rightarrow
ker\left(\pi_1A_2\rightarrow\pi_1M_1\right),  ker\left(\pi_1\partial A_1\rightarrow\pi_1M_1\right)\rightarrow
ker\left(\pi_1\partial A_2\rightarrow\pi_1M_1\right)$$ 
and $im\left(\pi_1A_1\rightarrow\pi_1M_1\right)$ is
amenable,
then $$\parallel M,\partial M\parallel\ge
\parallel M_1,\partial M_1\parallel.$$\end{lem}
\begin{pf} (i): For manifolds of dimensions $\le 2$ one
checks easily that there is no counterexample. So we are going
to assume that $n\ge3$.

We can restrict to the case that $A_1$ and $A_2$ are path-connected, since we may argue succesively
for their path-connected components. We denote by $A$ the image of $A_1$ in $M$.

The assumption $ker\left(\pi_1A_1\rightarrow\pi_1M_1\right)=
ker\left(\pi_1A_2\rightarrow\pi_1M_2\right)$ implies that $\pi_1M_1\rightarrow\pi_1M$ and $\pi_1M_2\rightarrow\pi_1M$ are injective. Thus $K\left(M_1\right),K\left(M_2\right)$ are submulticomplexes of $K\left(M\right)$.

For $i=1,2$ let $\partial_1M_i=M_i\cap\partial M$. (Thus $\partial M_i=\partial_1M_i\cup A_i$.) The assumption $ ker\left(\pi_1\partial A_1\rightarrow\pi_1M_1\right)
=ker\left(\pi_1\partial A_2\rightarrow\pi_1M_2\right)$
implies that $\pi_1\partial_1M_i\rightarrow\pi_1M$ is injective.

Let
$l_i:\left(M_i,\partial_1M_i\right)\rightarrow \left(M,\partial M\right)$ for $i=1,2$ be the inclusion.
 
We define bounded cohomology classes $\gamma_i\in H_b^n\left(M_i,\partial M_i\right)$ as follows. If $\parallel M_i,\partial M_i\parallel=0$, 
then we set $\gamma_i:=0$. If $\parallel M_i,\partial M_i\parallel\not=0$, 
then the cohomological fundamental class has a preimage $\beta_i\in H_b^n\left(M_i,\partial M_i\right)$ and we set
$\gamma_i:=\parallel M_i,\partial M_i\parallel \beta_i$.

In the main part of the proof we will now construct 
$\gamma\in H_b^n\left(M,\partial M\cup A\right)$ such that $\parallel\gamma
\parallel\le max\left\{\parallel\gamma_1\parallel,\parallel\gamma_2\parallel\right\}$ and $l_i^*\gamma=\gamma_i\in
H^n\left(M_i,\partial M_i\right)$.


We have the following sequence of commutative diagrams, where the first diagram is induced by inclusion and the second comes from naturality of $X\rightarrow K\left(X\right)$ and part ii) of \hyperref[iso3]{Proposition \ref*{iso3}}:
$$\begin{xy}\xymatrix{
\bigoplus_{i=1}^2 H_b^n\left(M_i,\partial M_i\right)\ar[d]^{incl^*\oplus incl^*} &H_b^n\left(M,\partial M\cup A\right)\ar[d]^{incl^*}\ar[l]^{l^*}\\
\bigoplus_{i=1}^2  H_b^n\left(M_i,\partial_1 M_i\right)\ar[d]^{J_1^*\oplus J_2^*} &H_b^n\left(M,\partial M\right)\ar[d]^{J^*}\ar[l]^{l_1^*\oplus l_2^*}\\
\bigoplus_{i=1}^2 H_b^n\left(K\left(M_i\right),K\left(\partial_1 M_i\right)\right)&H_b^n\left(K\left(M\right),K\left(\partial M\right)\right)\ar[l]^{k_1^*\oplus k_2^*}\\
\bigoplus_{i=1}^2 H^n_b\left(K\left(M_i\right),K\left(\partial_1 M_i\right)\right)^{\Pi_{M}\left(A\right)}
\ar[u]^{p_1\oplus p_2}&H_b^n\left(K\left(M\right),K\left(\partial M\right)\right)^{\Pi_M\left(A\right)}\ar[u]^{p}\ar[l]
}
\end{xy}$$


By part (i) of \hyperref[iso3]{Proposition \ref*{iso3}}, $J^*$ is an isometric isomorphism with inverse $I^*:H_b^n\left(K\left(M\right),K\left(\partial M\right)\right)
\rightarrow H_b^n\left(M,\partial M\right)$. 
By \hyperref[biggroup]{Lemma \ref*{biggroup}}, 
$\Pi_M\left(A\right)$ (as defined in Section 1.6.1) is amenable. Hence, by 
\hyperref[amen]{Lemma \ref*{amen}}(i), $p_1,p_2,p$
have left inverses $Av_1,Av_2,Av$ of norm $=1$.
For i=1,2 define $$\gamma^\prime_i:=Av_iJ_i^*incl^*\gamma_i
\in H_b^n\left(
K\left(M_i\right),
K\left(\partial_1 M_i\right)\right)^{\Pi_{M}\left(A\right)}.$$ They satisfy
$\parallel\gamma^\prime_i\parallel\le\parallel \gamma_i\parallel$.\\

In order to apply \hyperref[coramal]{Corollary \ref*{coramal}} to the pair $\left(K\left(M\right),K\left(\partial M\right)\right)$,
we have to check that the assumptions i)-v) of \hyperref[amalgamated]{Lemma \ref*{amalgamated}} are satisfied
with $$\left(K,K^\prime\right):=\left(K\left(M_1\right),K\left(\partial_1 M_1\right)\right), \left(L,L^\prime\right)=
\left(K\left(M_2\right),K\left(\partial_1 M_2\right)\right).$$

The only nonobvious conditions are possibly assumptions iv) and v). To check iv), let
$j:K\left(A\right)\rightarrow K\left(M\right)$ be the simplicial map defined in Section 1.3. Since $A$ is connected, we have that
$j\left(K\left(A\right)\right)
=K\left(\partial M_1\right)\cap K\left(\partial M_2\right)$ is a path-connected submulticomplex of $K\left(M\right)$. It is not hard to show, using the simplicial approximation theorem, that the homomorphisms
$\pi_1\left(j\left(K\left(A\right)\right)\right)\rightarrow  \pi_1\left(K\left(M_1\right)\right),
\pi_1\left(j\left(K\left(A\right)\right)\right)\rightarrow  \pi_1\left(K\left( M_2\right)\right)$ are injective. 

In particular, $\pi_1\left(K\left(M_1\right)\cup K\left(M_2\right)\right)$ is the amalgamated product of
$\pi_1\left(K\left(M_1\right)\right) $ and $\pi_1\left(K\left(M_2\right)\right)$, amalgamated over $\pi_1
j\left(K\left(A\right)\right)$. 

Again using the simplicial approximation theorem one can show that $\pi_1
j\left(K\left(A\right)\right)\simeq im\left(\pi_1A\rightarrow\pi_1M\right)$.
Thus $\pi_1\left(K\left(M_1\right)\cup K\left(M_2\right)\right)$ is isomorphic to
the amalgamated product of $\pi_1 M_1$ and $\pi_1 M_2$, amalgamated
over $im\left(\pi_1 A\right)$. By the Seifert-van Kampen Theorem the latter is isomorphic to $\pi_1M\simeq\pi_1K\left(M\right)$. This shows condition iv) and a similar 
argument proves condition v).

Thus we can apply \hyperref[coramal]{Corollary \ref*{coramal}} to $\gamma_1^\prime,\gamma_2^\prime$
and get $$\gamma^\prime\in H_b^n\left(K\left(M\right),
K\left(\partial M\right)\right)^{\Pi_M\left(A\right)},$$
satisfying $\parallel\gamma^\prime\parallel\le max\left\{
\parallel \gamma_1^\prime\parallel, \parallel \gamma_2^\prime\parallel
\right\}\le max\left\{\parallel \gamma_1\parallel,\parallel \gamma_2\parallel
\right\}$ and 
$k_i^*\gamma^\prime=\gamma_i^\prime$ for i=1,2.\\

If $c^\prime$ is a $G$-invariant relative cocycle representing $\gamma^\prime$, then $p^* c^\prime$ vanishes on simplices in $j\left(K\left(A\right)\right)$, by Observation 2. Since $I^*$ is induced by the chain homotopy equivalence $T:N\rightarrow K\left(N\right)$ and $T$ maps simplices in $A$ to simplices in $j\left(K\left(A\right)\right)$, we have that $I^*p^*c^\prime$ vanishes on simplices in $A$. That is, $I^*p^*c^\prime$ is a relative cocycle in $C_b^*\left(M,\partial M\cup A\right)$ and
$$I^*p^*\gamma^\prime      \in H_b^*\left(M,\partial M\right)$$ has a preimage
$$\gamma\in H_b^*\left(M,\partial M\cup A\right).$$
By construction it satisfies
$$\parallel\gamma\parallel\le \parallel \gamma^\prime\parallel\le
max\left\{\parallel\gamma_1\parallel,\parallel\gamma_2\parallel\right\}.$$

We want to check the equality $l_i^*\gamma=\gamma_i\in H^n\left(M_i,\partial M_i\right)$ for $i=1,2$. Since $dim\ H^n\left(M_i,\partial M_i\right)=1$, it suffices to 
evaluate both cohomology classes on a representative of $\left[M_i,\partial M_i\right]$. Choose some chain $z\in C_n\left(K\left(M_i\right)
\right)$ 
such that $\partial z=\partial_0z+\partial_1z$ with $\partial_0z$ representing $\left[A_i\right]\in H_{n-1}\left(A_i\right)$ and $\partial_1z$ representing $\left[\partial_1M_i
\right]\in H_{n-1}\left(\partial_1M_i\right)$. (In particular, $\partial z$ represents $\left[\partial M_i\right]\in H_{n-1}\left(M_i\right)$ and 
$z$ represents $\left[M_i,
\partial M_i\right]\in H_n\left(M_i,\partial M_i\right)$. Then, by \hyperref[obs2]{Observation \ref*{obs2}}, $\partial_0z\otimes 1=0$, thus
$z\otimes 1\in C_n\left(M_i,\partial_1M_i\right)\otimes_{{\Bbb Z}G}{\Bbb Z}$ 
is a relative cycle with $\gamma_i^\prime\left(z\otimes 1\right)=
\gamma_i\left(z\right)$. Since $r\left(l_i\left(z\right)\otimes 1\right)=z\otimes 1$, we compute
$$l_i^*\gamma\left(z\right)=\gamma\left(l_i\left(z\right)\right)=I^*p^*\gamma^\prime\left(l_i\left(z\right)\right)=\gamma^\prime\left(l_i\left(z\right)\otimes 1\right)=\gamma_i^\prime\left(z\otimes 1\right)=
p^*\gamma_i^\prime\left(z\right)=p^*Av\gamma_i\left(z\right).$$
 
It is easy to see that $\Pi_{M_i}\left(A_i\right)$ are connected and that the actions
of $\Pi_{M_i}\left(A_i\right)$ on $K\left(M_i\right)$ are continuous. 
Hence, all elements of $\Pi_{M_i}\left(A_i\right)$ are, as mappings from
$K\left(M_i\right)$ to itself, homotopic to the identity. From part (ii) of \hyperref[amen]{Lemma \ref*{amen}}, we get that $p^*\circ Av=id$ on $H_b^*\left(K\left(M_i\right),K\left(\partial_1M_i\right)\right)$.
Hence, we obtain $l_i^*\gamma=\gamma_i$ for $i=1,2$. 
 
Thus, we have constructed
$\gamma\in H_b^n\left(M,\partial M\cup A\right)$ such that $\parallel\gamma
\parallel\le max\left\{\parallel\gamma_1\parallel,\parallel\gamma_2\parallel\right\}$ and $l_i^*\gamma=\gamma_i\in
H^n\left(M_i,\partial M_i\right)$. We claim that this implies the wanted inequality. Indeed, in  $H_n\left(M,\partial M\cup A\right)$ we have $i_*\left[M,\partial M\right]={l_1}_*
\left[M_1,\partial M_1\right]+{l_2}_*\left[M_2,\partial M_2\right]$, which implies $i^*\gamma\left(\left[M,\partial M\right]\right)=
\gamma_1\left(\left[M_1,\partial M_1\right]\right)+
\gamma_2\left(\left[M_2,\partial M_2\right]\right)=
\sima+\simb.$ Thus, $\beta=\frac{1}{\sima+\simb}i^*\gamma$
is the relative fundamental cocycle of $\left(M,\partial M\right)$ and, by
$\parallel\beta\parallel\le\frac{1}{\sima +\simb}$ and duality follows
$\parallel M,\partial M\parallel\ge\parallel M_1,\partial M_1\parallel
+\parallel M_2,\partial M_2\parallel.$

(ii) In dimensions $\le 2$ we check
that the closed interval is the only connected counterexample. Assume then $n\ge3$.  Again we may suppose $A_1, A_2$ connected. We denote by $A$ the image of $A_1$ in $M$.  
Let $\partial_1M_1=M_1\cap\partial M$, which $\pi_1$-injects into $M_1$.
We consider the pair $$\left(N,N^\prime\right):=
\left(M_1\cup_{\left(id,0\right)+\left(f^{-1},1\right)}A_1\times\left[0,1\right],
\partial M_1\cup_{\left(id,0\right)+\left(f^{-1},1\right)}\partial\left(A_1\times\left[0,1\right]\right)\right)$$ and 
observe that $H_b^*\left(M,\partial M\cup A\right)=H_b^*\left(N,N^\prime\right)$.
Like in part (i), we get a commutative diagram
$$\begin{xy}\xymatrix{
H_b^n\left(M_1,\partial M_1\right)\ar[d]^{incl^*} &H_b^n\left(N,N^\prime\right)\ar[d]^{incl^*}\ar[l]^{l^*}\\
H_b^n\left(M_1,\partial_1 M_1\right)\ar[d]^{J_1^*} &H_b^n\left(N,\partial N\right)\ar[d]^{J^*}\ar[l]^{l^*}\\
H_b^n\left(K\left(M_1\right),K\left(\partial_1 M_1\right)\right)&H_b^n\left(K\left(N\right),K\left(\partial N\right)\right)\ar[l]^{k^*}\\
H^n_b\left(K\left(M_1\right),K\left(\partial_1 M_1\right)\right)^{\Pi_{M_1}\left(A_1\cup A_2\right)}
\ar[u]^{p_1}&H_b^n\left(K\left(N\right),K\left(\partial N\right)\right)^{\Pi_N\left(A\right)}\ar[u]^{p}\ar[l]
}
\end{xy}$$

Given $\gamma_1
\in H^*_b\left(M_1,\partial M_1\right)$,
define $$\gamma^\prime_1:=Av_1J_1^*incl^*\gamma_1
\in H_b^*\left(
(K\left(M_1\right)
,K\left(\partial M_1\right)\right)^{\Pi_{M_1}\left(A_1\cup A_2\right)}
.$$
Again, one checks easily that \hyperref[corhnn]{Corollary \ref*{corhnn}} can be applied to the pair 
$\left(K\left(N\right),K\left(\partial N\right)\right)$, namely that the assumptions of 
\hyperref[hnn]{Lemma \ref*{hnn}} are satisfied for $$\left(K,K^\prime\right):=\left(K
\left(M_1\right),K\left(\partial_1 M_1\right)\right), \left(L,L^\prime\right):=\left(h^{-1}j
\left(K\left(A_1\times\left[0,1\right]\right)\right),h^{-1}j\left(K\left(\partial A_1\times\left[0,1\right]\right)\right)\right).$$
From \hyperref[corhnn]{Corollary \ref*{corhnn}} we get $\gamma^\prime\in     H_b^n\left(K\left(N\right),K\left(\partial N\right)\right)^{\Pi_N\left(A\right)}.$
Define 
$\gamma:=I^*p\gamma^\prime\in H_b^n\left(N,N^\prime\right)=H_b^n\left(M,\partial M\cup A\right)$. Like in the proof of (i) we get that $j^*\gamma=\gamma_1$ and $\parallel \gamma\parallel\le\parallel \gamma_1\parallel$.
An argument analogous to the final step of i) shows that this implies $\parallel M,\partial M\parallel\ge\parallel M_1,\partial M_1\parallel$.\end{pf} 

\begin{lem}\label{Lemma9} 

i) Let $M_1,M_2$ be compact manifolds, $A_1$ resp.\ $A_2$ connected components
of $\partial M_1$ resp.\ $\partial M_2$ and assume that there exist connected sets
$A_i^\prime\subset M_i$ with $A_i^\prime\supset A_i$ and $\pi_1A_i^\prime$ amenable\footnotemark \footnotetext[12]{ The assumption of Lemma 8 is in particular satisfied if
$im\left(\pi_1\partial M_1\rightarrow\pi_1M_1\right)$ and
$im\left(\pi_1\partial M_2\rightarrow\pi_1M_1\right)$ are amenable and
the (singular) compression disks can be chosen to be disjoint. For
3-manifolds $M_i$, by a theorem of Jaco, cf.\ [4], there is $A_i^\prime\subset M_i$ with
$\pi_1A_i^\prime=im\left(\pi_1A_i\rightarrow\pi_1M_i\right)$
if $im\left(\pi_1A_i\rightarrow\pi_1M_i\right)$ is finitely presented.}.
Let $f:A_1\rightarrow A_2$ be a homeomorphism and
$M=M_1\cup M_2/f$ the glued manifold.
Then $$\parallel M,\partial M\parallel\le\parallel M_1,\partial M_1\parallel + \parallel M_2,\partial M_2\parallel.$$
ii) Let $M^\prime$ be a compact manifold, $A_1,A_2$ connected components of
$\partial M^\prime$ with $\pi_1A_i$ amenable. Let $f:A_1\rightarrow A_2$ be a homeomorphism and
$M=M^\prime/f$ the glued manifold. Then $$\parallel M,\partial M\parallel\le
\parallel M^\prime,\partial M^\prime\parallel.$$\end{lem}
\begin{pf} ii) is reduced to i) via the homeomorphism $M=M^\prime\cup_{\left(id,0\right)+\left(f,1\right)}\left(A_1\times I\right)$. (Note that $\parallel A_1\times I,A_1\times\left\{0,1\right\}\parallel =0$, since $\pi_1A_1$ is amenable.)
 
To prove i), we need the following reformulation of a theorem of Matsumoto-Morita.
For a space $X$  and $q\in N$ let $C_q\left(X\right)$ be the group of singular chains
and $B_q\left(X\right)$ the subgroup of boundaries. By Theorem 2.8.\ of \cite{mat} the following
two statements are equivalent:\\
a) there exists a number $K>0$ such that for any boundary $z\in B_q\left(X\right)$ there
is a
chain $c\in C_{q+1}\left(X\right)$ satisfying $\partial c=z$ and $\parallel c\parallel <
K\parallel z\parallel$,\\
b) the homomorphism $H_b^{q+1}\left(X\right)\rightarrow H^{q+1}\left(X\right)$ is injective.
 
Now let $\sum_{i=1}^ma_i\sigma_i$ and $\sum_{j=1}^nb_j\tau_j$
be representatives of $\left[M_1,\partial M_1\right]$ and
$\left[M_2,\partial M_2\right]$ with $$\sum_{i=1}^m\mid
a_i\mid\le\sima+\epsilon\mbox{\ \ \ and\ \ }\sum_{j=1}^n\mid
b_j\mid\le\simb+\epsilon.$$
By \hyperref[Prop4]{Proposition \ref*{Prop4}} we may suppose that
$$\parallel \partial\left( \sum_{i=1}^ma_i\sigma_i\right)\mid_{A_1}\parallel < \frac{\epsilon}{2K}\mbox{\ \ \ and\ \ }
\parallel \partial\left( \sum_{j=1}^nb_j\tau_j\right)\mid_{A_2}\parallel < \frac{\epsilon}{2K}.$$
(Indeed $\pi_1 A_i
\rightarrow\pi_1M_i$ factors over
$\pi_1A_i^\prime$, hence has amenable image.)
 
Let $A^\prime$ be the image of $A^\prime_1$ in $M$.
As $\pi_1A^\prime$ is amenable, $H_b^{q+1}\left(A^\prime\right)=0$ for $q\ge 0$ (\cite{gro},\cite{iva}),
hence, $H_b^{q+1}\left(A^\prime\right)\rightarrow H^{q+1}\left(A_i^\prime\right)$ is
injective and we get a constant $K$ with the property in a).

Therefore,
we find
$c\in C_*\left(A^\prime\right)\subset C_*\left(M\right)$ with $\parallel
c\parallel\le\epsilon$ and $$\partial c=\partial\left(\sum_{i=1}^ma_i\sigma_i+
\sum_{j=1}^nb_j\tau_j\right).$$
 
Then
$z=\sum_{i=1}^ma_i\sigma_i + \sum_{j=1}^nb_j\tau_j-c\in
C_*\left(\glu\right)$ is a fundamental
cycle with
$$ \parallel z\parallel < \parallel M_1,\partial M_1\parallel+\parallel M_2,\partial M_2\parallel+3\epsilon.$$
\end{pf}

\hyperref[lem 8]{Lemma \ref*{lem8}} and \hyperref[Lemma9]{Lemma \ref*{Lemma9}} imply {\bf Theorem 1} from the Introduction.

\begin{cor}\label{Cor4}

(i) Let A be a properly embedded
annulus in a compact 3-manifold M. If the identification $f:A_1\rightarrow A_2$ induces an isomorphism
$ker\left(\pi_1A_1\rightarrow\pi_1M_A\right)\approx
ker\left(\pi_1A_2\rightarrow\pi_1M_A\right)$ for the two 
images $A_1, A_2$ of $A$
in $M_A$, then\footnotemark \footnotetext[13]{If $\partial M$ consists of tori and A is an incompressible
annulus, then even $\parallel M_A,\partial M_A\parallel=\parallel
M,\partial M\parallel$ holds by a theorem of Soma (\cite{som1}).}
$$\parallel M_A,\partial M_A\parallel\le\parallel
M,\partial M\parallel.$$

(ii) Let T be an embedded torus in a compact 3-manifold M.
If the identification $f:T_1\rightarrow T_2$ induces an isomorphism
$ker\left(\pi_1T_1\rightarrow\pi_1M_T\right)\approx
ker\left(\pi_1T_2\rightarrow\pi_1M_T\right)$ for the two
images $T_1, T_2$ of $T$ in $M_T$, then 
$$\parallel M_T,\partial M_T\parallel
=\parallel M,\partial M\parallel.$$
\end{cor}

In \cite{som1}, a version of \hyperref[Cor4]{Corollary \ref*{Cor4}} has been proved for the special case that $\partial M$
consists of tori. The principal ingredient in the proof is the following statement:
\begin{pro}\label{Prop5} (\cite{som1}, Lemma 1):  Let $M$ be a compact 3-manifold whose boundary $\partial M$ consists of
tori and $H$ be a 3-dimensional compact submanifold of $int\left(M\right)$. Suppose
$int\left(H\right)$ is a hyperbolic 3-manifold and $\partial H$ is incompressible in $M$.Then we have $\parallel M,\partial M\parallel\ge\parallel H,\partial H\parallel$.\end{pro}

In \cite{thul}, this Proposition is stated for closed $M$ as Theorem 6.5.5. without writing a proof. In
\cite{som1}, it is then derived for
$M$ with toral boundary, using a doubling argument, see the proof of Lemma 1
in \cite{som1}. 
Hence, our proof seems to be the first written proof of \hyperref[Cor4]{Corollary \ref*{Cor4}} and \hyperref[Prop5]{Proposition \ref*{Prop5}}. (It
is easy to see that \hyperref[Cor4]{Corollary \ref*{Cor4}} implies \hyperref[Prop5]{Proposition \ref*{Prop5}}.)

We refer to \cite{som1} for the computation of the simplicial
volumes of all Haken 3-manifolds with (possibly empty) toral boundary, using Thurston's Theorem 6.5.5..
\subsection{Counterexamples}\label{examples}
The condition                            
$ker\left(\pi_1A_1\rightarrow\pi_1M_1\right)
\approx ker\left(\pi_1A_2\rightarrow\pi_1M_2\right)$ can not be weakened:
 
{\bf Example 1:} {\em Dehn fillings}\\
Let $K$ be a hyperbolic knot in $S^3$,
$V\subset S^3$ a regular neighborhood of $K$. Then
$$\parallel S^3\parallel=0<
\parallel S^3 - V,\partial V\parallel+\parallel V,\partial V\parallel,$$
although $\pi_1\partial V$ is amenable.
More generally, from Thurston's hyperbolic Dehn surgery theory we have $$\parallel M,\partial M\parallel < \parallel M - V,\partial V\cup\partial M\parallel=\parallel M - V,\partial V\cup\partial M\parallel+\parallel V,\partial V\parallel$$ holds if 
$M$ is obtained by performing hyperbolic Dehn filling at $M-V$.\\
\hyperref[lem8]{Lemma \ref*{lem8}} does not apply because the meridian of $\partial V$
is zero in $\pi_1V$, but not in
$\pi_1\left(M - V\right)$. \\
\\
The assumption '$A_i$ connected'
in \hyperref[Lemma9]{Lemma \ref*{Lemma9}}
can not be weakened:

{\bf Example 2:} {\em Heegaard splittings}\\
Let $M=H\cup H^\prime$ be a genus $g$ Heegaard splitting of a $3$-manifold $M$.
Let $D_1,\ldots,D_g\subset H^\prime$ be properly embedded disks with disjoint open regular neighborhoods $V_1,\ldots,V_g$ such that
$H^\prime - \cup_{i=1}^g V_i$ is a 3-ball $B$.
Let $A_i=V_i\cap\partial H^\prime$.

The manifold $M$ is then reconstructed as follows: $V_1,\ldots,V_g$ are glued to $H$ along the annuli $A_i$, afterwards
$B$ is glued along its whole boundary. 
Of course, $\parallel V_i,\partial V_i\parallel =0$ and $\parallel B,\partial B\parallel=0$.

If
\hyperref[Lemma9]{Lemma \ref*{Lemma9}} were applicable to the
annuli $A_i$, we would get that
$\parallel M,\partial M\parallel<\parallel H_g,\partial
H_g\parallel$. 

But there are $3$-manifolds of arbitrarily large
simplicial volume which admit Heegard splittings of a given genus.
To give an explicit example,
let $f$ be a pseudo-Anosov diffeomorphism on a surface of genus $g$,
and let $M_n$ be the mapping tori of the iterates $f^n$.
By Thurstons hyperbolization theorem, $M_1$ is hyperbolic. Hence,
$\parallel M_1\parallel>0$ and $\parallel M_n\parallel=n\parallel M_1\parallel$
becomes arbitrarily large. On the other hand, all $M_n$ admit a Heegaard splitting of
genus $2g+1$.

\subsection{Doubling 3-manifolds reduces simplicial volume}

For an oriented manifold $M$, let $DM$ denote the double of $M$, defined by glueing two
differently oriented copies of $M$ via the identity of $\partial M$. It is trivial that
$\parallel DM\parallel\le 2\parallel M,\partial M\parallel$. \hyperref[Theorem1]{Theorem \ref*{Theorem1}} implies,
of course:
if $M$ is a compact $3$-manifold with $\parallel \partial M\parallel=0$,
then $\parallel DM\parallel=2\parallel M,\partial M\parallel$. We will show that
this is actually an if-and-only-if condition.

\begin{thm}\label{Theorem2} If M is a compact $3$-manifold with incompressible boundary, then
$$\parallel\partial M\parallel >0 \Longleftrightarrow
\parallel DM\parallel < 2\parallel M,\partial M \parallel.$$\end{thm}
\begin{pf}
The claim can be reduced to compact irreducible manifolds because
$$\parallel M_1\sharp M_2,\partial\left(M_1\sharp M_2\right)\parallel
=\parallel M_1,\partial M_1\parallel + \parallel M_2,\partial M_2\parallel$$
holds
also for manifolds with boundary (of dimension $\ge3$). Indeed, defining a
fundamental class of the wedge $M_1\vee M_2$ as $\left[M_1\vee M_2,\partial
\left(M_1\vee M_2\right)\right]:={i_1}_*\left[M_1,\partial M_1\right]
+{i_2}_*\left[M_2,\partial M_2\right],$ for the inclusions $i_1:M_1\rightarrow
M_1\vee M_2$ and $i_2:M_2\rightarrow M_1\vee M_2$, we can define
the simplicial volume $\parallel M_1\vee M_2,\partial\left(M_1\vee M_2
\right)\parallel$ as the infimum over the $l^1$-norms of relative cycles
representing $\left[M_1\vee M_2,\partial\left(M_1\vee M_2\right)\right]$.
The proof of \hyperref[Theorem1]{Theorem \ref*{Theorem1}} shows that with this definition
$$\parallel M_1\vee M_2,\partial\left(M_1\vee M_2\right)\parallel=
\sima +\simb$$ holds.
Consider the projection from $M_1\sharp M_2$ to $M_1\vee M_2$ which 
pinches the connecting sphere to a point. It induces an isomorphism of
fundamental groups (this is the point, where dimension $\ge 3$ is needed)
and has degree 1. By the same argument as in the proof of \hyperref[degree]{Lemma \ref*{degree}},
we get $$\parallel M_1\sharp M_2,\partial\left(M_1\sharp M_2\right)\parallel
=\parallel M_1\vee M_2,\partial\left(M_1\vee M_2\right)\parallel,$$
which proves the first equality.
An obvious generalization of this argument shows 
$$\parallel D\left(M_1\sharp M_2\right)\parallel = \parallel DM_1\parallel
+ \parallel DM_2\parallel.$$

Any compact $3$-manifold with incompressible boundary is a connected sum of $3$-manifolds which are either irreducible
or $S^1\times S^2$. The above equalities (and $\parallel S^1\times S^2\parallel=0$) reduce the claim to irreducible $3$-manifolds.
 
So assume $M$ irreducible. By standard $3$-manifold theory
we can cut $M$ along properly embedded disks and incompressible 
properly embedded annuli and tori such that the obtained pieces are either 
Seifert fibered, I-bundles or "simple" in the sense of \cite{neus}. 

We argue that these pieces $M_i$ satisfy the conclusion of Theorem 2. For a Seifert
fibration, the boundary consists of tori, hence, there is nothing to prove.
If $M_i$ is an $I$-bundle, then $DM$ is an $S^1$-bundle and $\parallel DM_i\parallel=0$ holds by Corollary 6.5.3 of \cite{thul}. 

Finally, if 
$M_i$ is "simple" and $\parallel\partial M_i\parallel >0$, then $M_i$ admits a hyperbolic metric with totally geodesic boundary (and possibly cusps), and
the totally geodesic boundary is non-empty. 
If there are no cusps, then $\parallel DM_i\parallel < 2\parallel M_i,\partial M_i\parallel$ follows from \cite{jun}.
If there are cusps, then the same inequality follows from Theorem 6.3.\ in \cite{k}.

If $M$ is a compact irreducible $3$-manifold with $\parallel \partial M\parallel>0$, then clearly $\parallel \partial M_i\parallel > 0$ holds for at least one of
the pieces in its Jaco-Shalen-Johannson decomposition. To conclude the
proof
of Theorem 2, we still need  to prove that the implication
$$\parallel D\left(M_F\right)\parallel < 2\parallel M_F,\partial M_F\parallel \Longrightarrow
\parallel DM\parallel < 2\parallel M,\partial M\parallel$$
(with 
$M_F$ defined as in the Introduction) holds true whenever
$M$ is a compact $3$-manifold with incompressible boundary
and $F$ an incompressible, properly embedded annulus, torus or disk. 

If $F$ is a torus, this claim follows easily from \hyperref[Cor4]{Corollary \ref*{Cor4}}(ii). 

If $F=A$ is an annulus, the
claim follows from the
somewhat paradoxical observation that 
$$\parallel M_A,\partial M_A\parallel
\le\parallel M,\partial M\parallel,\mbox{\ but\ }\parallel D\left(M_A\right)\parallel
\ge\parallel DM\parallel.$$
 
Indeed, $DM_A$ is obtained from $DM$ by cutting off
an incompressible torus
$T=DA$ and identifying afterwards in a different way the
pairs of incompressible annuli which are "halves" of
the same copy of T. In other words $$\left(D\left(M_A\right)\right)_{2A}=\left(DM\right)_T.$$
Hence, we get $$\parallel D\left(M_A\right)\parallel \ge \parallel\left(DM\right)_T\parallel
=\parallel DM\parallel,$$
what implies the claim.
If $F$ is a disk, the argument is the same.\end{pf}\\
\\
{\bf Handlebodies}. \\
Note that for any $n$-dimensional compact manifold $M$
$$\parallel M,\partial M\parallel\ge\frac{1}{n}\parallel \partial M\parallel$$ holds. Namely, the 
boundary operator $$\partial:H_n\left(M,\partial M\right)\rightarrow H_{n-1}\left(
\partial M\right)$$ maps the relative fundamental class $\left[M,\partial M\right]$
to the fundamental class $\left[\partial M\right]$. It is obvious that $$\parallel\partial\parallel\le n+1.$$ 
(Each $n$-simplex has $n+1$ faces.) For a representative $\sum_{i=1}^ra_i\sigma_i$ of $\left[M,\partial M\right]$                                          
one even gets $$\parallel \partial\sum_{i=1}^ra_i\sigma_i\parallel\le n\parallel\sum_{i=1}^r
a_i\sigma_i\parallel,$$
because each $\sigma_i$ has to have at least one face not in $\partial M$,
cancelling against some other face. (In \cite[Proposition 2.7]{bfp} this is improved to $\parallel\partial M\parallel\le (n-1)\parallel M,\partial M\parallel$.)

Let $H_g$ denote the $3$-dimensional handlebody of genus $g$. $H_g$ is a (g-1)-fold
covering of $H_2$, hence 
$$\parallel H_g,\partial H_g\parallel =C\left(g-1\right).$$
From the above argument follows $C\ge\frac{4}{3}$. (In fact, we showed in \cite{kue}, by constructing
triangulations of handlebodies, that $\frac{4}{3}\le C\le 3$. According to the recent preprint \cite{bfp} the accurate value is $C=3$.)
 
The double of $H_g$ is the $g$-fold connected sum $S^2\times S^1\sharp\ldots\sharp S^2\times S^1,$
whose simplicial volume vanishes.
Thus there are manifolds
$M$ of arbitrarily large simplicial volume with $\parallel DM\parallel =0$.\\
\\
In general, a compact $3$-manifold $M$ can be cut along disks into finitely many pieces 
$M_j$ which have incompressible boundary, and it is not hard to show that
$\parallel DM\parallel=0$
if and only if all these $M_j$ have a
Jaco-Shalen-Johannson-decomposition without simple pieces in the sense
of \cite{neus}.

\noindent
Thilo Kuessner,
Korea Institute for Advanced Study,
85 Hoegi-ro, Dongdaemun-gu,
130-722 Seoul,
Republic of Korea


\begin{thebibliography}{29}
\bibitem{bbfipp}M.\ BUCHER, M.\ BURGER, R.\ FRIGERIO, A.\ IOZZI, C.\ PAGLIANTINI and M.\ B.\ POZZETTI, 'Isometric embeddings in bounded cohomology.' {\em J. Topol. Anal.} 6  (2014) 1-25.
\bibitem{bfp}M.\ BUCHER, R.\ FRIGERIO and C.\ PAGLIANTINI, 'The simplicial volume of 3-manifolds with boundary';  http://arxiv.org/abs/1208.0545
\bibitem{gro}M.\ GROMOV, 'Volume and Bounded
Cohomology', {\em Public.\ Math.\ IHES} 56 (1982) 5-100.
\bibitem{iva}N.\ IVANOV, 'Foundations of the theory of
bounded cohomology', {\em J.\ Sov.\ Math.}\ 37 (1987) 1090-1114.
\bibitem{jun}D.\ JUNGREIS, 'Chains that realize the Gromov
invariant of hyperbolic manifolds',  {\em Ergodic Theory and Dynamical Systems} 17 (1997) 643-648.
\bibitem{kk}S.\ KIM and T. KUESSNER, 'Simplicial volume of compact manifolds with amenable boundary'. J.\ Topol.\ Anal.\ 7 (2015), no. 1, 23-46.
\bibitem{kue}T.\ KUESSNER, 'Gromov volume of compact manifolds', Diplomarbeit, FU Berlin, 1996
\bibitem{k}T.\ KUESSNER, 'Efficient fundamental cycles of
cusped hyperbolic manifolds', {\em Pac.\ J.\ Math.}\ 211 (2003) 283-314.
\bibitem{k4}T.\ KUESSNER, 'Generalizations of Agol's inequality and nonexistence of tight laminations', {\em Pac.\ J.\ Math.}\ 251 (2011) 109-172.
\bibitem{lam}K.\ LAMOTKE, {\em Semisimpliziale algebraische Topologie}, (Die Grundlehren der mathematischen Wissenschaften, Springer-Verlag, Berlin, 1968).
\bibitem{mat}S.\ MATSUMOTO and S.\ MORITA, 'Bounded Cohomology of
Certain Groups of Homeomorphisms', {\em Proc.\ Amer.\ Math.\ Soc.}\ 94 (1985) 539-544.
\bibitem{may}J.\ P.\ MAY, {\em Simplicial objects in algebraic topology} (Chicago Lect.\ Math., UCP, 1992).
\bibitem{neus}W.\ NEUMANN and G.\ SWARUP, 'Canonical decompositions of
3-manifolds', {\em Geometry and Topology} 1 (1997) 21-40.
\bibitem{park}H.\ PARK, 'Relative bounded cohomology', {\em Top.\ Appl.}\ 131 (2003) 203-234.
\bibitem{som1}T.\ SOMA, 'The Gromov invariant of links', 
{\em Invent.\ Math.}\ 64 (1981) 445-454.
\bibitem{thul}W.\ THURSTON, 'The Geometry and Topology of
3-Manifolds', Lecture Notes; http://msri.org/publications/books/gt3m 
\end{thebibliography}
\end{document}